\documentclass[preprint, 11pt]{elsarticle}

\usepackage[T1]{fontenc}
\usepackage[utf8]{inputenc}
\usepackage[hmargin=0.9in]{geometry}
\usepackage[babel,german=quotes]{csquotes}
\usepackage{subcaption}
\usepackage{graphicx}
\usepackage[colorlinks=true,citecolor=blue,urlcolor=blue]{hyperref}
\usepackage{caption}
\usepackage{amsmath}
\usepackage{amsthm}
\usepackage{amssymb}
\usepackage{bm, amsfonts}
\usepackage[dvipsnames]{xcolor}
\usepackage{soul,framed}
\usepackage{tikz}
\usetikzlibrary{plotmarks}
\usepackage{bm}
\usepackage{makecell}
\usepackage{titlesec}
\usepackage{multirow}
\usepackage{indentfirst}
\usepackage{aligned-overset}

\newtheorem{theorem}{Theorem}[section]

\theoremstyle{definition}

\theoremstyle{remark}
\newtheorem{remark}[theorem]{Remark}
\newtheorem*{theorem*}{Remark}

\numberwithin{equation}{section}

\usepackage[notref,notcite,final]{showkeys} 
\usepackage[mode=multiuser]{fixme} 
\FXRegisterAuthor{nm}{envar}{NM}
\FXRegisterAuthor{dk}{envdk}{DK}
\FXRegisterAuthor{jv}{envjv}{JV}
\FXRegisterAuthor{jm}{envjm}{JM}
\fxusetheme{color}
\colorlet{fxwarning}{blue}

\makeatletter \def\pgfuseplotmark#1{\pgftransformresetnontranslations\csname pgf@plot@mark@#1\endcsname} \makeatother

\newcommand{\vertiii}[1]{{\left\vert\kern-0.25ex\left\vert\kern-0.25ex\left\vert #1 
		\right\vert\kern-0.25ex\right\vert\kern-0.25ex\right\vert}}

\colorlet{shadecolor}{green}

\titleformat{\paragraph}
{\normalfont\normalsize\itshape}{\theparagraph}{1em}{}
\titlespacing*{\paragraph}
{0pt}{3.25ex plus 1ex minus .2ex}{1.5ex plus .2ex}

\makeatletter
\def\ps@pprintTitle{%
	\let\@oddhead\@empty
	\let\@evenhead\@empty
	\def\@oddfoot{
		\footnotesize\itshape
		\ifx\@journal\@empty Elsevier
		\else\@journal\fi
		\hfill\today
	}%
	\let\@evenfoot\@oddfoot}
\makeatother

\begin{document}
	\begin{frontmatter}
		\title{A structure-preserving Lagrangian discontinuous Galerkin method using flux and slope limiting}
  		\author{Joshua Vedral$^{1}$}
		\ead{joshua.vedral@math.tu-dortmund.de}
		\author{Nathaniel Morgan$^{2}$\corref{cor1}}
		\ead{nmorgan@lanl.gov}
		\cortext[cor1]{Corresponding author}
		\author{Dmitri Kuzmin$^{1}$}
		\ead{kuzmin@math.uni-dortmund.de}
  		\author{Jacob Moore$^{3}$}
		\ead{jlmoore@iser.msstate.edu}
		
		\address{$^{1}$Institute of Applied Mathematics (LS III), TU Dortmund University\\ Vogelpothsweg 87,
			D-44227 Dortmund, Germany}
        \address{$^{2}$Computational Physics Division, Los Alamos National Laboratory\\ Bikini Atoll Rd SM 30, Los Alamos, NM 87545, USA}
		\address{$^{3}$ISER Mississippi State University \\ 3909 Halls Ferry Rd, Vicksburg, MS 39180, USA}

		\journal{Preprint submitted to Elsevier}
		
		\begin{abstract}
		We introduce a Lagrangian nodal discontinuous Galerkin (DG) cell-centered hydrodynamics method for solving multi-dimensional hyperbolic systems. By incorporating an adaptation of Zalesak's flux-corrected transport algorithm, we combine a first-order positivity-preserving scheme with a higher-order target discretization. This results in a flux-corrected Lagrangian DG scheme that ensures both global positivity preservation and second-order accuracy for the cell averages of specific volume. The correction factors for flux limiting are derived from specific volume and applied to all components of the solution vector. We algebraically evolve the volumes of mesh cells using a discrete version of the geometric conservation law (GCL). The application of a limiter to the GCL fluxes is equivalent to moving the mesh using limited nodal velocities.  Additionally, we equip our method with a locally bound-preserving slope limiter to effectively suppress spurious oscillations. Nodal velocity and external forces are computed using a multidirectional approximate Riemann solver to maintain conservation of momentum and total energy in vertex neighborhoods. Employing linear finite elements and a second-order accurate time integrator guarantees GCL consistency. The results for standard test problems demonstrate the stability and superb shock-capturing capabilities of our scheme.
		\end{abstract}
		
		
	\end{frontmatter}
	\section{Introduction}
	Numerical solutions to conservation laws in gas dynamics and magneto-hydrodynamics can be obtained using finite difference (FD), finite volume (FV), or finite element (FE) methods that can be further subdivided into Eulerian, Lagrangian, and arbitrary Lagrangian-Eulerian (ALE) approaches. Within the Eulerian framework, the mesh remains fixed and the material moves through the mesh, introducing numerical viscosity. Especially in the presence of shocks and contact discontinuities, this behavior has a negative impact on accuracy. In contrast, Lagrangian methods directly associate the material motion with the computational mesh. In these methods, the mesh moves alongside the material, offering a flexible framework for capturing shocks in multidimensional flows.
	
	Within the Lagrangian framework, two distinct concepts have emerged --- staggered grid and cell-centered schemes. The first Lagrangian schemes were one-dimensional FD staggered-grid hydrodynamic (SGH) methods that discretized the conservation equations in a spatially staggered manner \cite{morgan2021}.  In this approach, the momentum is computed at the nodes of the mesh, while thermodynamic variables are located at the cell center. In their seminal work, von Neumann and Richtmyer introduced a FD SGH method \cite{morgan2021,vonneumann1950} for one-dimensional flows and equipped it with artificial viscosity, which was later extended by Kolsky to the two-dimensional case \cite{morgan2021}.  Many subsequent works in the field, including those by Wilkins \cite{wilkins1963} and Burton \cite{Burton1994}, would advance the FV Lagrangian SGH method, which uses staggered control volumes.  Across multi-dimensional Lagrangian SGH methods, artificial viscosity is essential to achieving stable and accurate solutions of shock problems. Thus, many research efforts have focused on developing accurate and robust artificial viscosity-based methods, wherein dissipation is either added explicitly \cite{kuropatenko1967,morgan2013,wilkins1980} or through solving Riemann problems \cite{christensen1990,dukowicz1985,loubere2010,loubere2010a,loubere2013,maire2011,morgan2014}. 
	
	In contrast to FV SGH methods, the cell-centered hydrodynamics (CCH) approach solves the conservation equations for each cell within a single control volume. The computation of nodal velocity for mesh movement is achieved by solving multidimensional approximate Riemann problems at each node or by least squares fitting the velocities on the cell faces that come from 1D Riemann solvers. The first FV CCH method was proposed by Godunov \cite{godunov1999, godunov1979}. In contrast to the SGH framework, approximate Riemann solvers of CCH schemes inherently introduce sufficient numerical dissipation. Thus, there is no need for adding artificial viscosity terms. Research in early to mid 1980's by Addessio et al. \cite{addessio1986} yielded a multi-dimensional Lagrangian CCH method based on using 1D Riemann solvers and calculating the nodal mesh velocity by least-squares fitting the 1D face velocities.  This least squares fitting approach could give rise to spurious mesh motion, requiring the use of a vorticity correction technique for stable solutions \cite{dukowicz1992}.  Loub\'ere et al. \cite{loubere2004} introduced a novel cell-centered Lagrangian scheme based on a fully Lagrangian formulation of the gas dynamics equations, i.e., the gradient and divergence operators are expressed in Lagrangian form. Despr\'es and Mazeran \cite{despres2005} extended this scheme, resulting in the pivotal GLACE scheme \cite{carre2009}. They assumed the velocity to be continuous at nodes, allowing them to develop the first node-centered approximate Riemann solver. The discrete volume evolution equation thus satisfies the geometric conservation law (GCL). The resulting two-dimensional Lagrangian FV CCH scheme is conservative and entropy consistent. Despite its strengths, the nodal velocity in GLACE schemes is sensitive to cell aspect ratios, leading to numerical instabilities. To address this, a new nodal solver was developed by introducing multiple pressure fluxes at the grid nodes, giving rise to the EUCCLHYD scheme proposed by Maire et al. \cite{maire2007,maire2008}. Several robust nodal Riemann solvers can be found, e.g., in \cite{burton2013,maire2009,maire2007,morgan2013a}.  Nodal multidimensional approximate Riemann solvers are the mainstay of Lagrangian CCH methods as they deliver accurate and stable mesh motion. 
	
	The blending of finite element methods with Lagrangian schemes dates back to the two-dimensional Lagrangian DG hydrodynamics (DGH) scheme which was proposed by Loub\'ere et al. \cite{loubere2004} for linear triangular meshes. While their scheme uses nodal basis functions, nearly all subsequent efforts focused on modal DGH techniques; see, e.g., \cite{jia2011,lieberman2019,liu2018,liu2019,morgan2018,vilar2012a,vilar2014}.  Recently, an arbitrary-order nodal Lagrangian DGH scheme was developed by Moore et al. \cite{moore2021} for smooth flows and weak shocks.  This 3D nodal Lagrangian DGH scheme was later combined with a classical slope limiter to simulate shock problems with linear hexahedral elements \cite{MultiMat2022Presentation}. However, aggressive limiting was required to keep the internal energy positive, degrading the accuracy of simulations.  The modal and nodal Lagrangian DGH methods can be viewed as extensions of the Lagrangian FV CCH methods.  
    Ellis \cite{ellis2010} 
    proposed a natural extension of Lagrangian SGH approaches by employing biquadratic polynomial basis functions for kinematic variables and bilinear basis functions for thermodynamic variables. Dobrev et al. \cite{dobrev2011} extended the methodology employed by Ellis to allow for arbitrary order of kinematic and thermodynamic basis functions \cite{dobrev2012}, as well as for axisymmetric problems \cite{dobrev2013}. Scovazzi et al. \cite{scovazzi2007,scovazzi2010} proposed a SUPG-stabilized formulation for Lagrangian hydrodynamics providing a globally conservative approximation. Their approach allows for the representation of pressure gradients on element interiors, unlike traditional CCH schemes. 
	
	Inspired by the progress on 2D Lagrangian DGH schemes, and seeing literature gaps on robust 3D Lagrangian DGH methods, we propose a new element-centered Lagrangian DGH method for solving the Euler equations of gas dynamics in 2D and 3D Cartesian coordinates.  Following the principles of flux-corrected transport (FCT) methods \cite{boris1973,kuzmin2012,kuzmin2002,lohner1987,zalesak1979}, we combine a first-order positivity-preserving scheme with a high-order counterpart. Employing flux limiting guarantees global positivity preservation for the element (also called cell or zone) averages of specific volume, while slope limiting effectively suppresses spurious oscillations. Our proposed scheme is fully conservative, GCL consistent, and second-order accurate. 
	
	We begin, in Section \ref{sec:goveqs}, by formulating the governing equations of gas dynamics in Lagrangian form. The DG spatial discretization is presented in Section \ref{sec:dg}. We proceed, in Section \ref{sec:rie}, with the description of a multidimensional approximate Riemann solver. Our limiting strategy is presented in Section \ref{sec:lim}. A brief discussion on GCL consistency follows in Section \ref{sec:gcl}. The physical quantities are evolved in time using the time integrator presented in Section \ref{sec:time}. Numerical experiments for standard benchmark problems are performed in Section \ref{sec:numex}. Finally, we give concluding remarks in Section \ref{sec:concl}.
	\section{Governing equations}
	\label{sec:goveqs}
	The differential Lagrangian forms of the evolution equations for specific volume, momentum, and specific total energy of an inviscid compressible fluid are given by
	\begin{align}
	\begin{split}
	\varrho\frac{\mathrm{d}\nu}{\mathrm{d}t}&=\nabla \cdot \mathbf{u},\\
	\varrho\frac{\mathrm{d}\mathbf{u}}{\mathrm{d}t}&=\nabla \cdot \boldsymbol{\sigma}, \\
	\varrho\frac{\mathrm{d}\tau}{\mathrm{d}t}&=\nabla \cdot (\boldsymbol{\sigma} \cdot \mathbf{u}),
	\end{split}
	\label{eq:lagr}
	\end{align}
	where $\varrho$ is the density, $\nu$ is the specific volume, $\mathbf{u}$ is the velocity, $\tau$ is the specific total energy, and $\boldsymbol{\sigma}$ is the stress. In gas dynamics, the stress tensor is isotropic and reads $\boldsymbol{\sigma}=-p\mathbf{I}$, where $p$ denotes the pressure and $\mathbf{I}$ is the identity matrix. The pressure is calculated using an equation of state for the material. In the case of a polytropic ideal gas, the pressure can be calculated using 
	\begin{align}
		p=(\gamma-1)\varrho e,
		\label{eq:eos}
	\end{align}
	where $\gamma$ is the heat capacity ratio and $e$ is the specific internal energy. 
	
	We further denote by $\mathbf{x}$ the position of a particle in the Eulerian frame and by $\mathbf{\tilde{x}}$ the position of a particle in the Lagrangian frame. The trajectory of a particle or a mesh node moving with velocity $\mathbf{u}$ can now be defined as the solution of the initial value problem
	\begin{align}
	\frac{\mathrm{d}\mathbf{x}}{\mathrm{d}t}=\mathbf{u}, \qquad \mathbf{x}|_{t=0}=\mathbf{\tilde{x}}.
	\label{eq:link}
	\end{align}
	\section{Discontinuous Galerkin discretization}
	\label{sec:dg}
	Let $\{\omega_c\}$, $c=1,\ldots,C_h$ be a decomposition of the solution domain $\omega(t)$ into $C_h$ non-overlapping polygonal elements (also called cells). We denote by $m_c$ the constant mass in the element $\omega_c$. The initial element configuration, denoted by $\omega_c^0$, $c=1,\ldots,C_h$ satisfies $\cup_{c\in C_h}\omega_c^0=\omega(0)=:\omega^0$. Each cell has $N_p$ nodes (or vertices). Let $N_c$ represent the number of unknowns in a cell $\omega_c$. The boundary of the element $\omega_c$ is denoted by $\partial \omega_c$.
	
	The deformation of each element is described by linking Lagrangian coordinates (the initial configuration) with Eulerian coordinates (the configuration at time $t$). Using \eqref{eq:link}, the mapping from the initial configuration to the current configuration can be written as
	\begin{align*}
        \mathbf{x}=\Pi(\mathbf{\tilde{x}},t),
	\end{align*}
	where $\Pi:\omega_c^0\to  \omega_c$ is the mapping function. Similarly, an isoparametric mapping from a reference element $\omega_R$ to the current configuration is
	\begin{align*}
	  \mathbf{x}=\Phi(\boldsymbol{\xi},t)
          =\sum_{p=1}^{N_p}b_p(\boldsymbol{\xi})\mathbf{x}_p(t),\qquad \boldsymbol{\xi}\in\omega_R.
	\end{align*}
	Here, $\Phi:\omega_{R}\to \omega_c$ is the mapping function, $\mathbf{\xi}$ are the reference coordinates, $b_p$ are canonical shape functions at the nodes, and $\mathbf{x}_p$ are the nodal coordinates. Both the deformation gradient $\mathbf{F}=\frac{\partial \mathbf{x}}{\partial \mathbf{\tilde x}}$ and the Jacobian matrix $\mathbf{J}=\frac{\partial \mathbf{x}}{\partial \mathbf{\boldsymbol{\xi}}}$ are time-dependent.
	
	
	Unlike existing DGH methods in which specific volume, velocity and specific total energy are represented as modal fields approximated with Taylor polynomials \cite{lieberman2019,liu2019,morgan2018}, we opt for nodal basis functions to be determined later. In the process of numerical integration, we evaluate the conserved quantities at Gaussian quadrature points. 
	
	We write the evolution equations \eqref{eq:lagr} in the generic form
	\begin{align}
	\varrho \frac{\mathrm{d}\mathbb{U}}{\mathrm{d}t}=-\nabla \cdot \mathbb{H}(\mathbb{U}),
	\label{eq:evol}
	\end{align}
	where $\mathbb{U}=(\nu,\mathbf{u},\tau)^{\top}$ is the vector of conserved quantities and $\mathbb{H}=(-\mathbf{u},-\boldsymbol{\sigma},-\boldsymbol{\sigma}\cdot \mathbf{u})$ is the corresponding flux vector.
	
	Proceeding with spatial discretization using DG finite elements, we introduce the vector of conserved unknowns $\mathbb{U}_h=(\nu_h,\mathbf{u}_h,\tau_h)^{\top}$. On each element $\omega_c$, we seek an approximate solution 
	\begin{align}
		\mathbb{U}_{h,c}(\mathbf{x},t):=\mathbb{U}_h(\mathbf{x},t)|_{\omega_c}=\sum_{j=1}^{N_c}\mathbb{U}_j(t)\psi_j(\mathbf{x}), \quad \mathbf{x}\in \omega_c, \quad t\in[0,T],
		\label{eq:sol}
	\end{align}  
	where $\psi_j$, $j=1,\ldots,N_c$ are linear Lagrange basis functions.
	
	By inserting \eqref{eq:sol} into \eqref{eq:evol}, then multiplying by a sufficiently smooth test function $\Psi_i$ and integrating over the element $w_c$, we derive the local weak formulation
	\begin{align*}
	  \int_{\omega_c}\Psi_i \Bigg(\varrho \frac{\mathrm{d}\mathbb{U}_{h,c}}{\mathrm{d}t}+\nabla \cdot \mathbb{H}(\mathbb{U}_{h,c})\Bigg)\,\mathrm{d}w=
          \int_{\partial \omega_c} \Psi_i[\mathbb{H}(\mathbb{U}_{h,c})
            \cdot\mathbf n-
            \mathbb{H}^\ast(\mathbf{u}^\ast,\boldsymbol{\sigma}^\ast;\mathbf{n})]
            \,\mathrm{ds},
	\end{align*}
	where $\mathbf{n}$ is the outward unit normal to the surface and $\mathbb{H}^\ast$ is a numerical flux that depends on an approximate solution of a Riemann problem posed on the element boundary $\mathrm{ds}$. The computation of the Riemann velocity $\mathbf{u}^\ast$ and stress $\boldsymbol{\sigma}^\ast$ will be discussed in the next section. 
        
	Using the local consistent mass matrix 
	\begin{align}
	  M_c=(m_{ij})_{i,j=1,\ldots,N_c},\qquad m_{ij}=
          \int_{\omega_c} \varrho \Psi_i\Psi_j \,\mathrm{d}w,
		\label{eq:mass}
	\end{align}
	we obtain
	\begin{align*}
	  \sum_{j=1}^{N_c}m_{ij}\frac{\mathrm{d}\mathbb{U}_j}{\mathrm{d}t}=
          -\int_{\omega_c}\Psi_i\nabla \cdot \mathbb{H}(\mathbb{U}_{h,c})\,\mathrm{d}w+
\int_{\partial \omega_c} \Psi_i[\mathbb{H}(\mathbb{U}_{h,c})
            \cdot\mathbf n-
            \mathbb{H}^\ast(\mathbf{u}^\ast,\boldsymbol{\sigma}^\ast;\mathbf{n})]
            \,\mathrm{ds}, \qquad i=1,\ldots,N_c.
	\end{align*}
 	Employing integration by parts and the divergence theorem yields
 	\begin{align}
 		\sum_{j=1}^{N_c}m_{ij}\frac{\mathrm{d}\mathbb{U}_j}{\mathrm{d}t}=\int_{\omega_c}\nabla \Psi_i\cdot \mathbb{H}(\mathbb{U}_{h,c})\,\mathrm{d}w-\int_{\partial \omega_c} \Psi_i \mathbb{H}^\ast(\mathbf{u}^\ast,\boldsymbol{\sigma}^\ast; \mathbf{n})\,\mathrm{ds}, \qquad i=1,\ldots,N_c.
 		\label{eq:dg}
 	\end{align}

	To compute the integrals that appear in \eqref{eq:dg} on the reference element $\omega_R$, we express the corresponding evolution equations as
	\begin{align}
		\sum_{j=1}^{N_c}m_{ij}\frac{\mathrm{d}\mathbb{U}_j}{\mathrm{d}t}=\int_{\omega_R}\nabla_{\xi} \Psi_i\cdot \mathbb{H}(\mathbb{U}_{h,c})\cdot \mathbf{J}^{-1}\det(\mathbf{J})\,\mathrm{d}w_\mathrm{R}-\int_{\partial \omega_R} \Psi_i \cdot\mathbb{H}^\ast(\mathbf{u}^\ast,\boldsymbol{\sigma}^\ast;\mathbf{n})\cdot \mathbf{J}^{-1}\det(\mathbf{J}) \,\mathrm{ds_R}
        \label{eq:dgref}
	\end{align}
	for $i=1,\ldots,N_c$, where $\mathrm{d}w_\mathrm{R}$ and $\mathrm{ds_R}$ denote integration over the reference element and its surface, respectively.
	
	It is straightforward to verify that the nodal Lagrangian DG method \eqref{eq:dg} conserves mass, momentum, and total energy. Summing over the basis functions and using their partition of unity property, we obtain
	\begin{align}
		\sum_{i=1}^{N_c}\sum_{j=1}^{N_c}m_{ij} \frac{\mathrm{d}\mathbb{U}_j}{\mathrm{d}t}=-\int_{\partial \omega_c} \mathbb{H}^\ast(\mathbf{u}^\ast,\boldsymbol{\sigma}^\ast;\mathbf{n}) \,\mathrm{ds}.\label{eq:cellsum}
	\end{align}
	Summing over all elements and using the definition of the mass matrix in \eqref{eq:mass}, we have
	\begin{align*}
		\sum_{c=1}^{C_h}\int_{\omega_c}\varrho \frac{\mathrm{d}\mathbb{U}_{h,c}}{\mathrm{d}t}\,\mathrm{d}w=0.
	\end{align*}
	\begin{remark}
		The element averages of the Lagrangian DG method \eqref{eq:dg} correspond to a first-order FV method (cf. \cite{gallice2022,georges2016}).  
		\label{rem:avgs}
	\end{remark}
	\section{Riemann problem}
	\label{sec:rie}
\begin{figure}
        \centering
        \includegraphics[trim={0 2.9in 0 2.9in},clip,width=0.5\linewidth]{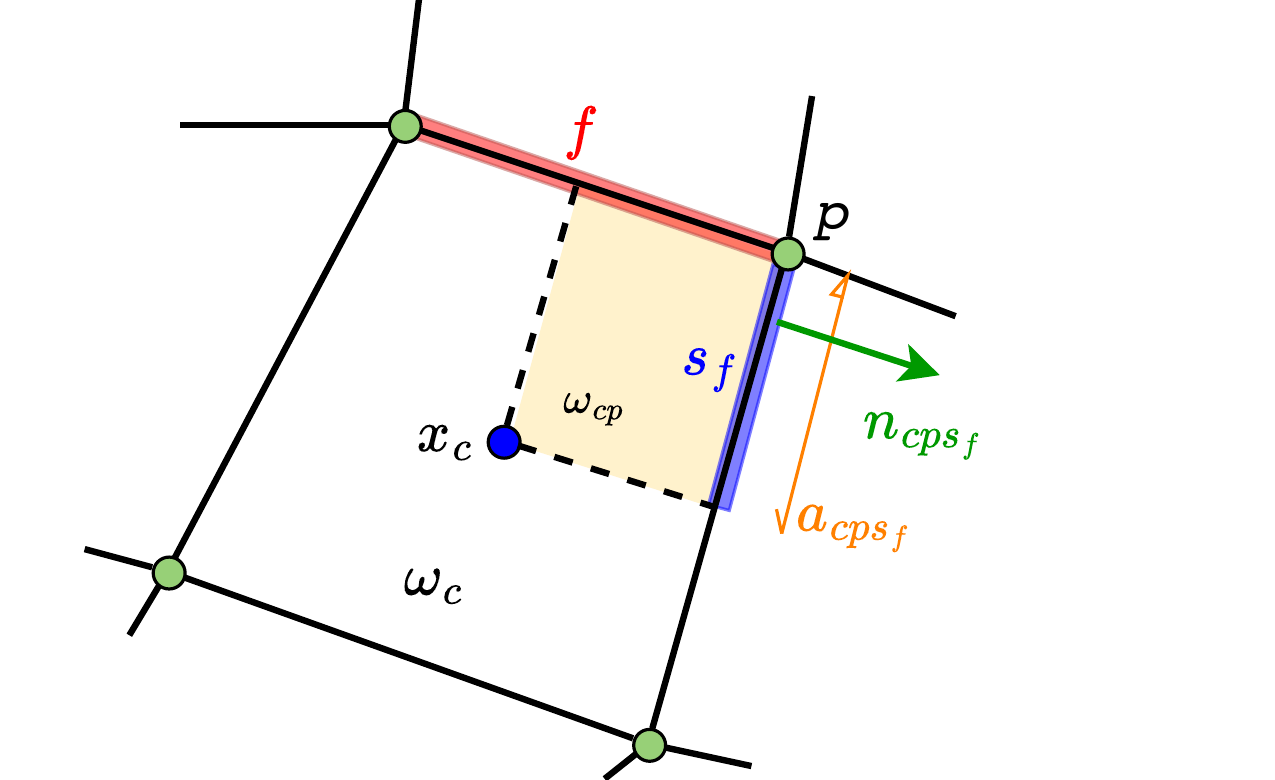}
        \caption{Geometrical entities attached to element $\omega_c$.}
        \label{fig:geom}
    \end{figure}

	Before starting with the description of the Riemann solver, we introduce some notation. The element $\omega_c$ is defined by its nodes (or vertices) $p\in \mathcal{P}(c)$, where $\mathcal{P}(c)$ is the set of all nodes of $\omega_c$. We denote by $\mathcal{C}(p):=\{\omega_c: \; p \in \omega_c\}$ the set of elements sharing the node $p$. We split $\omega_c$ into subcells $\omega_{cp}$ (cf. \cite{gallice2022,georges2016,liu2019,morgan2018}). Here, $\omega_{cp}$ is the subcell associated with cell $\omega_c$ and node $p$, formed by joining the element centroid $\mathbf{x}_c$ to the midpoints of faces $f$ attached to node $p$. Hence, the set of subcells $\omega_{cp}$ for $p\in \mathcal{P}(c)$ constitutes a partition of the element $\omega_c$, i.e., $\omega_c=\cup_{p\in \mathcal{P}(c)}$$\omega_{cp}$. Let $p$ be a node of element $\omega_c$ and $f$ be a face connected to node $p$ in cell $\omega_c$. We denote by $\mathcal{F}(c,p)$ the set of faces of element $\omega_c$ attached to $p$. Each face $f$ of element $\omega_c$ is divided into subfaces $s_f$ via the partition of $\omega_c$ induced by the subcells $\omega_{cp}$ for $p\in \mathcal{P}(c)$. We denote by $\mathcal{SF}(c,p)$ the set of subfaces of element $\omega_c$ attached to corner $l_{cp}$, equivalent to the set of faces of subcell $\omega_{cp}$ connected to node $p$. We then define $a_{cps_f}$ as the weighted length of the subface $s_f$ related to node $p$ in element $\omega_c$, and $a_{cps_f}\mathbf{n}_{cps_f}$ as the weighted area normal vector to the subface $s_f$ related to node $p$ in element $\omega_c$. The geometrical entities attached to element $\omega_c$ are displayed in Fig. \ref{fig:geom}.

	We employ the nodal multidirectional approximate Riemann solver proposed in \cite{morgan2018} to solve the Riemann problem arising in \eqref{eq:dg}. This solver, based on work by Burton et al. \cite{burton2013}, computes the Riemann force acting on a subface of the element by
	\begin{align}
		\mathbf{F}_{cps_f}^\ast=a_{cps_f}\mathbf{n}_{cps_f}\boldsymbol{\sigma}_{cps_f}^\ast=a_{cps_f}\mathbf{n}_{cps_f} \boldsymbol{\sigma}_{l_{cp}}+\mu_{cps_f}a_{cps_f}(\mathbf{u}_p^\ast-\mathbf{u}_{l_{cp}}).
		\label{eq:rieforce}
	\end{align}
	Here, $\mu_{cps_f}=\varrho_{cps_f} c^+ |\hat{\bf e}_{cp} \cdot {\bf n}_{cps_f}|$, where $\varrho_{cps_f} c^+$ denotes the acoustic impedance, $c^+$ is the sound speed of the element, and $\hat{\bf e}_{cp}$ is a unit vector in the direction of the difference between the node velocity and the element average velocity.

	Momentum conservation at node $p$ requires the summation of all forces around the node to be zero, i.e., 
	\begin{align}
		\sum_{\omega_c\in \mathcal{C}(p)}\sum_{s_f\in \mathcal{SF}(c,p)}\mathbf{F}_{cps_f}^\ast=0.
		\label{eq:momcons}
	\end{align}
	Rearranging \eqref{eq:rieforce} and using \eqref{eq:momcons} yields the Riemann velocity
	\begin{align}
		\mathbf{u}_p^\ast=\sum_{\omega_c\in \mathcal{C}(p)}\sum_{s_f\in \mathcal{SF}(c,p)}\frac{\mu_{cps_f} a_{cps_f}\mathbf{u}_{l_{cp}}-a_{cps_f}\mathbf{n}_{cps_f}\boldsymbol{\sigma}_{l_{cp}}}{\mu_{cps_f}a_{cps_f}}.
		\label{eq:rievel}
	\end{align}
	The Riemann stress $\boldsymbol{\sigma}_{cps_f}^\ast$ can now be calculated using \eqref{eq:rievel} and the Riemann jump condition. 
	
	The approximate Riemann problem in \eqref{eq:rieforce} introduces dissipation based on discontinuities in the velocity and/or stress fields. Combining \eqref{eq:momcons} and \eqref{eq:rievel}, we observe conservation of total energy. Liu et al. \cite{liu2018} verified that the Lagrangian DG method along with this Riemann solver satisfies the second law of thermodynamics.
	\section{Limiting}
	\label{sec:lim}
	The Lagrangian DG discretization \eqref{eq:dg} combined with the Riemann solver \eqref{eq:rieforce} suffers from oscillations near shocks and steep gradients. To address these issues, we present a novel limiting approach for Lagrangian DGH methods. Using the limiting framework outlined in \cite{kuzmin2021}, we incorporate flux correction for element averages and slope correction for directional derivatives. Aligning with the principles of flux-corrected transport (FCT) methods \cite{boris1973,kuzmin2012,kuzmin2002,lohner1987,zalesak1979}, we adapt Zalesak's algorithm \cite{zalesak1979,zalesak2005} to ensure non-negativity of element-averaged specific volumes in our DG discretization. The correction factors derived from specific volume are then applied to all components of the solution vector. To ensure local bound preservation, we make use of a clip-and-scale limiter \cite{anderson2017,kuzmin2023,lohmann2017}. In line with the employed flux limiting approach, the correction factors for slope limiting are also computed from the specific volume and applied to all fields.
	\subsection{Element averages of a Lagrangian DGH method}
	As stated in Remark \ref{rem:avgs}, our Lagrangian DG method reduces to a first-order FV scheme when considering only the element averages. Let  $\bar{\mathbb{U}}_{h,c}=\frac{1}{m_c}\int_{\omega_c}\mathbb{U}_h\,\mathrm{d}w$ denote the cell average of $\mathbb{U}_h$ over $\omega_c$. In view of \eqref{eq:cellsum}, the
        element averages of a Lagrangian DG approximation satisfy \cite{gallice2022,georges2016}
	\begin{align}
	  m_c\frac{\mathrm{d}\bar{\mathbb{U}}_{h,c}}{\mathrm{d}t}+\int_{\partial \omega_c}\mathbb{H}^\ast(\mathbf{u}^\ast,\boldsymbol{\sigma}^\ast;\mathbf{n})
          \,\mathrm{ds}=0.
		\label{eq:avgssd}
	\end{align}
 Employing a first-order explicit time integrator turns \eqref{eq:avgssd} into
	\begin{align}
	  m_c(\bar{\mathbb{U}}_{h,c}^{n+1}-\bar{\mathbb{U}}_{h,c}^n)+\Delta t\int_{\partial \omega_c}\mathbb{H}^\ast(\mathbf{u}^\ast,\boldsymbol{\sigma}^\ast;\mathbf{n}) \,\mathrm{ds}=0,
		\label{eq:avgsd}
	\end{align}
	where $\bar{\mathbb{U}}_{h,c}^{n+1}$ and $\bar{\mathbb{U}}_{h,c}^n$ denote the approximation of $\bar{\mathbb{U}}_{h,c}(t)$ at time instants $t^{n+1}$ and $t^n$, respectively, and $\Delta t$ is the time step size, i.e., $\Delta t=t^{n+1}-t^n$.
	
	Applying the subcell-based partition of $\omega_c$ from \cite{gallice2022,georges2016} and adopting the notation from the previous section, we rewrite the surface integral in \eqref{eq:avgsd} as
	\begin{align}
	  \int_{\partial \omega_c}
\mathbb{H}^\ast(\mathbf{u}^\ast,\boldsymbol{\sigma}^\ast;\mathbf{n})\,\mathrm{ds}&=\sum_{p\in \mathcal{P}(c)}\sum_{s_f\in \mathcal{SF}(c,p)}a_{cps_f}\mathbb{H}_{cps_f}\cdot \mathbf{n}_{cps_f}=0,
		\label{eq:avgssub}
	\end{align}
	where $\mathbb{H}_{cps_f}=(-\mathbf{u}_p,p_{cps_f}\mathbf{I},p_{cps_f}\mathbf{u}_p)^\top$. Inserting \eqref{eq:avgssub} into \eqref{eq:avgsd}, the element averages of a Lagrangian DG approximation can be evolved by
	\begin{align}
		m_c(\bar{\mathbb{U}}_{h,c}^{n+1}-\bar{\mathbb{U}}_{h,c}^n)+\Delta t\sum_{p\in \mathcal{P}(c)}\sum_{s_f\in \mathcal{SF}(c,p)}a_{cps_f}\mathbb{H}_{cps_f}\cdot \mathbf{n}_{cps_f}=0.
		\label{eq:avgsst}
	\end{align}

	\subsection{Flux-corrected transport}
	In FCT methods, a property-preserving low-order scheme is combined with a high-order target discretization, which generally violates physical properties. Following this approach, a low-order version of \eqref{eq:avgsst} can be derived using
	\begin{align}
		\mathbf{u}_p^L=\mathbf{u}_p^n, \qquad \mathbf{n}_{cps_f}^L=\mathbf{n}_{cps_f}^n.
		\label{eq:fctlo}
	\end{align}
	The corresponding low-order fluxes are denoted by $\mathbb{H}_{cps_f}^L$. The so-defined approximation $\bar{\mathbb{U}}_h^L$ is positivity preserving \cite{georges2016} and first-order accurate in time. Alternatively, a second-order accurate but generally not positivity-preserving approximation $\bar{\mathbb{U}}_h^H$ can be obtained using the temporal average
	\begin{align}
		\mathbf{u}_p^H=\frac{\mathbf{u}_p^{n+1}+\mathbf{u}_p^n}{2}, \qquad \mathbf{n}_{cps_f}^H=\frac{\mathbf{n}_{cps_f}^{n+1}+\mathbf{n}_{cps_f}^n}{2}.
		\label{eq:fctho}
	\end{align}
	The corresponding high-order fluxes are denoted by $\mathbb{H}_{cps_f}^H$.
	
	The low- and high-order approximations derived from \eqref{eq:fctlo} and \eqref{eq:fctho}, respectively, are linked by 
	\begin{align*}
		m_c\bar{\mathbb{U}}_{h,c}^H=m_c\bar{\mathbb{U}}_{h,c}^L+\Delta t \sum_{p\in \mathcal{P}(c)}\mathbb{H}_{cps_f}^B,
	\end{align*}
	where
	\begin{align*}
		\mathbb{H}_{cps_f}^B=\sum_{s_f\in \mathcal{SF}(c,p)}a_{cps_f}(\mathbb{H}_{cps_f}^L\cdot \mathbf{n}_{cps_f}^L-\mathbb{H}_{cps_f}^H\cdot \mathbf{n}_{cps_f}^H)
	\end{align*}
	are numerical fluxes that possess the conservation property 
	\begin{align*}
		\sum_{\omega_c\in \mathcal{C}(p)}\mathbb{H}_{cps_f}^B=0.
	\end{align*}
	Now, we can formulate a flux-corrected Lagrangian DG scheme (cf. \cite{kuzmin2021})
	\begin{align}
		m_c\bar{\mathbb{U}}_{h,c}^{n+1}=m_c\bar{\mathbb{U}}_{h,c}^L+\Delta t \sum_{p\in \mathcal{P}(c)}\alpha_p\mathbb{H}_{cps_f}^B.
		\label{eq:fct}
	\end{align}
	The conservation property is preserved by defining
	\begin{align*}
		\alpha_p=\min_{\omega_c\in \mathcal{C}(p)} \alpha_{cp}
	\end{align*}
	using correction factors $\alpha_{cp}\in [0,1]$ that ensure global positivity preservation for the element averages. 
	\subsection{Flux limiting}
	It remains to show how the correction factors $\alpha_{cp}$ are calculated. In essence, any element-based limiting technique outlined in \cite[Chapter 4]{kuzmin2023} can be adapted and applied. However, instead of constraining element contributions to a node, we constrain nodal contributions to an element. Scalar limiting is applied to the specific volume, and all fluxes in \eqref{eq:fct} are scaled by the resulting correction factor.
	
	Let $\mathbb{H}_{cps_f}^B|_{1/\varrho}:=\mathbb{H}_{cps_f}^B(\mathbf{u}^B,\cdot,\cdot)^\top$ be the flux function entries corresponding to the specific volume evolution equation. The nodal-based multidimensional version of Zalesak's limiting strategy \cite{kuzmin2023,zalesak1979,zalesak2005} is as follows:\\
	(1) Compute the sums of antidiffusive fluxes:
	\begin{align}
		P_c^+=\sum_{p\in \mathcal{P}(c)}\max(0,\mathbb{H}_{cps_f}^B|_{1/\varrho}), \qquad P_c^-=\sum_{p\in \mathcal{P}(c)}\min(0,\mathbb{H}_{cps_f}^B|_{1/\varrho}).
		\label{eq:zal1}
	\end{align}
	(2) Determine the bounds for limited increments:
	\begin{align}
		Q_c^+=\gamma_c(\mathbb{U}_c^{\max}|_{1/\varrho}-\bar{\mathbb{U}}_{h,c}^L|_{1/\varrho}), \qquad Q_c^-=\gamma_c(\mathbb{U}_c^{\min}|_{1/\varrho}-\bar{\mathbb{U}}_{h,c}^L|_{1/\varrho}).
		\label{eq:zal2}
	\end{align}
	(3) Compute the bounds for correction factors:
	\begin{align}
		R_c^+ = \min\Bigg(1,\frac{Q_c^+}{P_c^+}\Bigg), \qquad R_c^-=\min\Bigg(1,\frac{Q_c^-}{P_c^-}\Bigg).
		\label{eq:zal3}
	\end{align}
	(4) Select an upper bound:
	\begin{align}
		\alpha_{cp}=\begin{cases}
		R_c^+ \; &\text{if} \; \mathbb{H}_{cps_f}^B|_{1/\varrho}\geq 0, \\
		R_c^- \; &\text{if} \; \mathbb{H}_{cps_f}^B|_{1/\varrho}< 0.
		\label{eq:zal4}
		\end{cases}
	\end{align}
	 We set the upper and lower \textit{global} bounds for the specific volume as $\mathbb{U}_c^{\max}|_{1/\varrho}:=\mathbb{U}_c^{\max}(1/\varrho,\cdot,\cdot)^\top=\infty$ and $\mathbb{U}_c^{\min}|_{1/\varrho}:=\mathbb{U}_c^{\min}(1/\varrho,\cdot,\cdot)^\top=0$, respectively.
	Additionally, we set $\gamma_c=m_c/\Delta t$. The correction factors computed from \eqref{eq:zal1}--\eqref{eq:zal4} ensure a positivity-preserving approximation for the element averages of specific volume.
	
	\begin{remark}
		Applying a flux limiter is equivalent to moving the node $p$ using the velocity
		\begin{align}
		\mathbf{u}_{p,\mathrm{lim}}=\alpha_p\mathbf{u}_p^H+(1-\alpha_p)\mathbf{u}_p^L.
		\label{eq:meshvel}
		\end{align}
	\end{remark}
	\begin{remark}
		Our approach involves scalar limiting of the specific volume. Alternatively, one can employ FCT-type sequential limiting \cite[Chapter 4]{kuzmin2023}, where different correction factors are computed for each field. For the Euler equations of gas dynamics, after the density limiting step, the gradients of total energy and momentum are adjusted to account for the density changes \cite{dobrev2018}. An extension of this method to the equations of ideal magnetohydrodynamics can be found in \cite{kuzmin2020c}. 
	\end{remark}
	\subsection{Slope limiting}
	The flux limiter ensures the non-negativity of the element-averaged specific volume $\bar{\mathbb{U}}_{h,c}|_{1/\varrho}$, $c=1,\ldots,C_h$. However, local spurious oscillations, while not violating the positivity of the solution, may still arise. To eliminate these oscillations, we construct a locally bound-preserving conservative approximation $\mathbb{U}_{h,c}^*$ using the clip-and-scale limiting strategy \cite{anderson2017,kuzmin2023,lohmann2017}.
	
	The element averages $\bar{\mathbb{U}}_{h,c}$ of a Lagrangian DG approximation can be expressed in terms of their nodal values as follows:
	\begin{align*}
		\bar{\mathbb{U}}_{h,c}=\frac{1}{m_c}\sum_{j=1}^{N_c}m_{c,j}\mathbb{U}_j, \qquad 	m_{c,j}=\int_{\omega_c}\varrho_{h}\Psi_j\,\mathrm{d}w,
	\end{align*}
	where $ \sum_{j=1}^{N_c}m_{c,j}=m_c$. By defining $\mathbb{H}_{c,j}:=m_{c,j}(\mathbb{U}_j-\bar{\mathbb{U}}_{h,c})$, the nodal values $\mathbb{U}_j$ are given by
	\begin{align*}
		\mathbb{U}_j=\bar{\mathbb{U}}_{h,c}+\frac{\mathbb{H}_{c,j}}{m_{c,j}}, \qquad \sum_{j=1}^{N_c}\mathbb{H}_{c,j}=0.
	\end{align*}
	We aim to construct a locally bound-preserving conservative approximation
	\begin{align}
		\mathbb{U}_j^\ast=\bar{\mathbb{U}}_{h,c}+\frac{\mathbb{H}_{c,j}^\ast}{m_{c,j}}, \qquad \sum_{j=1}^{N_c}\mathbb{H}_{c,j}^\ast=0
		\label{eq:slapprox}
	\end{align}
	using limited fluxes $\mathbb{H}_{c,j}^\ast$ to be defined below.
	
	Again, we denote by $\mathbb{H}_{c,j}|_{1/\varrho}:=\mathbb{H}_{c,j}(\mathbf{u},\cdot,\cdot)^\top$ the flux function entries corresponding to the specific volume evolution equation. Following the flux limiting approach, the correction factors for slope limiting are derived from the specific volume and subsequently applied to all components of the solution vector.
	
	We begin by limiting each component of the element flux individually, i.e., 
	\begin{align*}
		\tilde{\mathbb{H}}_{c,j}|_{1/\varrho}=\max(\mathbb{H}_{c,j}^{\min}|_{1/\varrho},\min(\mathbb{H}_{c,j}|_{1/\varrho},\mathbb{H}_{c,j}^{\max}|_{1/\varrho})),
	\end{align*}
	where 
	\begin{align}
		\mathbb{H}_{c,j}^{\min}|_{1/\varrho}:=\gamma_{c,j}(\mathbb{U}_{c,j}^{\min}|_{1/\varrho}-\mathbb{U}_j|_{1/\varrho}), \qquad \mathbb{H}_{c,j}^{\max}|_{1/\varrho}:=\gamma_{c,j}(\mathbb{U}_{c,j}^{\max}|_{1/\varrho}-\mathbb{U}_j|_{1/\varrho})
		\label{eq:locbounds}
	\end{align}
	and $\gamma_{c,j}=m_{c,j}/\Delta t$. 
	Next, we compute the sums of clipped element contributions
	\begin{align}
		P_c^+=\sum_{j=1}^{N_c}\max(0,\tilde{\mathbb{H}}_{c,j}|_{1/\varrho}), \qquad P_c^-=\sum_{j=1}^{N_c}\min(0,\tilde{\mathbb{H}}_{c,j}|_{1/\varrho})
		\label{eq:cas1}
	\end{align}
	and subsequently scale positive or negative components as follows:
	\begin{align}
		\mathbb{H}_{c,j}^\ast|_{1/\varrho}= \begin{cases}
		-\frac{P_c^-}{P_c^+}\tilde{\mathbb{H}}_{c,j}|_{1/\varrho} \quad &\text{if} \; \tilde{\mathbb{H}}_{c,j}|_{1/\varrho}>0 \; \text{and} \; P_c^+ + P_c^->0, \\
		-\frac{P_c^+}{P_c^-}\tilde{\mathbb{H}}_{c,j}|_{1/\varrho} \quad &\text{if} \; \tilde{\mathbb{H}}_{c,j}|_{1/\varrho}<0 \; \text{and} \; P_c^+ + P_c^-<0, \\
		\phantom{-\frac{P_c^+}{P_c^-}}\tilde{\mathbb{H}}_{c,j}|_{1/\varrho} \quad &\text{otherwise}.
		\end{cases}
		\label{eq:cas2}
	\end{align}
	We set the upper and lower \textit{local} bounds for the specific volume as
	\begin{align*}
		\mathbb{U}_{c,j}^{\max}|_{1/\varrho}:=\mathbb{U}_{c,j}^{\max}(1/\varrho,\cdot,\cdot)^\top=\max_{\sum_{i\in \mathcal{N}_j}}\mathbb{U}_i|_{1/\varrho}, \qquad \mathbb{U}_{c,j}^{\min}|_{1/\varrho}:=\mathbb{U}_{c,j}^{\min}(1/\varrho,\cdot,\cdot)^\top=\min_{\sum_{i\in \mathcal{N}_j}}\mathbb{U}_i|_{1/\varrho},
	\end{align*}
	respectively. Here, $\mathcal{N}_j$ is the index set of nodes belonging to elements containing node $j$.

	\begin{remark}
		The limiting step \eqref{eq:cas2} enforces the zero-sum condition in \eqref{eq:slapprox}.
	\end{remark}
	
	\begin{remark}
		The clip-and-scale strategy is first applied to the specific volume, with the resulting correction factors then used across all fields. Alternatively, an extended version of the clip-and-scale limiter designed for systems can be employed. In this approach, after the scalar limiting of the specific volume, a product rule variant of the clip-and-scale limiter is applied to the remaining fields. For further details, we refer the reader to \cite[Section 6.3.2]{kuzmin2023}.
	\end{remark}
	\section{GCL consistency}
	\label{sec:gcl}
	As explained in \cite{loubere2008}, the volume $v_c^{n+1}$ of an element can be calculated using its definition as a function of vertex coordinates or evolved using the geometric conservation law (GCL). Solving \eqref{eq:dg} involves a discretized GCL. The two approaches are equivalent for certain types of elements if the mesh velocity is defined appropriately. To maintain volumetric consistency, we evolve the nodes according to \eqref{eq:meshvel}.
	
	In the one-dimensional case, the discretized GCL for an element $\omega_c=[x_i,x_{i+1}]$ reads
	\begin{align*}
		\tilde{v}_{i+1/2}^{n+1}=\tilde{v}_{i+1/2}^n+\Delta t(u_{i+1}-u_i),
	\end{align*}
	whereas the geometric definition of the evolving element volume is
	\begin{align*}
		v_{i+1/2}(t)=x_{i+1}(t)-x_i(t), \qquad t \in [t^n,t^{n+1}].
	\end{align*}
	Let the nodal velocities be fixed on $[t^n,t^{n+1}]$. Then we have
	\begin{align*}
		x_i^{n+1}=x_i^n+\Delta t u_i, \qquad x_{i+1}^{n+1}=x_{i+1}^n+\Delta t u_{i+1}
	\end{align*}
	and therefore
	\begin{align*}
		v_{i+1/2}^{n+1}=v_{i+1/2}^n+\Delta t(u_{i+1}-u_i).
	\end{align*}
	Thus, the discretized GCL holds for geometric volumes.
	
	To show GCL consistency for triangles and quadrilaterals, we use the fact that \cite[Section 1.3.2]{georges2016}
	\begin{align}
		\frac{\mathrm{d}v_c}{\mathrm{d}t}=\sum_{p\in \mathcal{P}(c)}\sum_{s_f\in \mathcal{SF}(c,p)}a_{cps_f}\mathbf{u}_{p}\cdot\mathbf{n}_{cps_f}.
		\label{eq:gcltriag}
	\end{align}
	Since the coordinates $x_p(t)$ and $y_p(t)$ of node $p$ are linear in $t$, the time derivative of $v_c(t)$ and hence the right-hand side of \eqref{eq:gcltriag} is also linear in $t$. Therefore, any second-order accurate time discretization ensures GCL consistency. However, a first-order accurate time integrator may not guarantee it. A more sophisticated proof of GCL consistency for triangles and quadrilaterals can be found in the \hyperref[sec:appendix]{Appendix}.
	
	\begin{remark}
		The impact of the time discretization on GCL consistency and stability is further discussed in \cite{farhat2001, nobile1999}. Increasing the degree of the finite element mapping results in a higher polynomial degree of $v_c(t)$. Consequently, the order of the time integrator must match the order of the space discretization for a Lagrangian DG method to be GCL consistent.
	\end{remark}

    \begin{remark}
        GCL consistency in 3D has been investigated by Georges \cite{georges2016} and Georges et al. \cite{georges2016a} who proposed face-splitting techniques applicable to arbitrary polyhedra. Polyhedra, by definition, have planar faces. In contrast, general hexahedral elements defined by trilinear mappings typically have non-planar faces. The assumption of planar faces simplifies the analysis, as it enables consistent and exact computation of the element volume using methods such as tetrahedral decomposition or closed-form geometric formulas. For elements with non-planar faces, the geometry is not equivalent to a collection of tetrahedra, and the geometric volume defined by the finite element mapping may differ from that obtained by a piecewise linear decomposition.
    \end{remark}
	\section{Time integration}
    \label{sec:time}
    The Lagrangian DG scheme \eqref{eq:dg} can be written in matrix form as
    \begin{align}
        M\frac{\mathrm{d}\mathbb{U}_h}{\mathrm{d}t}=R(\mathbb{U}_h),
    \end{align}
    where $M$ denotes the global consistent mass matrix and $R(\mathbb{U}_h)$ is the residual vector. 

    The numerical solution is evolved in time using the second-order strong stability preserving (SSP) Runge-Kutta method \cite{gottlieb2001}
    \begin{align*}
        \mathbb{U}_h^{(1)}&=\mathbb{U}_h^n+\Delta tM^{-1}R(\mathbb{U}_h^n),\\
        \mathbb{U}_h^{n+1}&=\frac{1}{2}\mathbb{U}_h^n+\frac{1}{2}\mathbb{U}_h^{(1)}+\frac{1}{2}\Delta tM^{-1}R(\mathbb{U}_h^{(1)}).
    \end{align*}
    The nodal positions $\mathbf{x}_p$ are updated using
    \begin{align*}
        \mathbf{x}_p^{(1)}&=\mathbf{x}_p^n+\Delta t\mathbf{u}_{p,\mathrm{lim}}^n,\\
        \mathbf{x}_p^{n+1}&=\frac{1}{2}\mathbf{x}_p^n+\frac{1}{2}\mathbf{x}_p^{(1)}+\frac{1}{2}\Delta t\mathbf{u}_{p,\mathrm{lim}}^{(1)},
    \end{align*}
    where $\mathbf{u}_{p,\mathrm{lim}}$ denotes the limited Riemann velocity at node $p$, see \eqref{eq:meshvel}.

    We refer the reader to \cite{georges2016,maire2009} for how to evaluate the time step size based on the CFL criterion.
	\section{Numerical examples}
	\label{sec:numex}
	We solve a series of standard test problems for the Euler equations of gas dynamics to demonstrate the accuracy and robustness of our DGH method. Flux and slope limiters are applied at the end of each Runge-Kutta stage, and all test problems use the gamma-law equation of state as defined in \eqref{eq:eos}.  A novel contribution of this work is the development of a Lagrangian DGH method capable of simulating three-dimensional problems involving both smooth and shock-driven flows. To the best of our knowledge, accurate and robust 3D  Lagrangian DGH schemes for shock-dominated problems are lacking in the literature. 
	
	For smooth flows, we consider the Taylor-Green vortex problem to verify the second-order accuracy of the fully discrete scheme. For shock-driven flows, we study the Sedov blast wave, the Noh problem, and the triple point problem, which together demonstrate the accuracy and stability of our DGH method in the presence of strong shocks. 
	
	In some test cases, we introduce an additional elementwise-constant parameter $\beta_c$ that is applied after slope limiting to prevent spurious mesh motion and mesh tangling. Specifically, after obtaining the slope-limited solution $\mathbb{U}_j^*$ from \eqref{eq:slapprox}, we further modify it as
	\begin{align*}
		\mathbb{U}_j^{*,new} = \beta_c \mathbb{U}_j^* \qquad \forall j \in N_c.
	\end{align*}
	By choosing $\beta_c=1$, we recover exactly the slope-limited solution $\mathbb{U}_j^*$ from \eqref{eq:slapprox}. A similar idea was investigated by Liu et al. \cite{liu2018}, who adjusted the slope-limiting parameter in a Barth-Jespersen limiter and found that the expected order of convergence is preserved for $\beta_c\in[0.6,1.0]$. Following their findings, we also present results using $\beta_c=0.8$ in cases where mesh tangling is observed. We emphasize that the scheme remains IDP regardless of the value of $\beta_c$; the parameter is introduced solely to enhance mesh stability. 
	
	\subsection{2D Taylor-Green vortex problem}
	In our first numerical experiment, we investigate the two-dimensional Taylor-Green vortex problem \cite{dobrev2012,morgan2015,taylor1937}, which can be derived by adapting an analytical solution of the incompressible Navier-Stokes equations to the compressible Euler equations. The computational domain is $\Omega=(0,1)^2$, with $\gamma=5/3$. The initial conditions are given by
	\begin{align*}
		\varrho^0 = 1, \quad \mathbf{u}^0 = [\sin(\pi x)\cos(\pi y), -\cos(\pi x )\sin(\pi y)]^{\top}, \quad p^0 = \frac{1}{4}[\cos(2\pi x)+\cos(2\pi y )] + 1.
	\end{align*}
	To satisfy the total energy evolution equation, an energy source term $S_E$ is added:
	\begin{align*}
		S_E = \frac{\pi}{4}\frac{(\varrho^0)^3}{\gamma-1}[\cos(3\pi x)\cos(\pi y)-\cos(\pi x)\cos(3\pi y)].
	\end{align*} 
	For details on the derivation and construction of this source term, we refer the reader to \cite{vilar2012a}. 
	
\begin{figure}
		\centering
		
		\begin{subfigure}[b]{0.19\textwidth}
			\includegraphics[width=\linewidth]{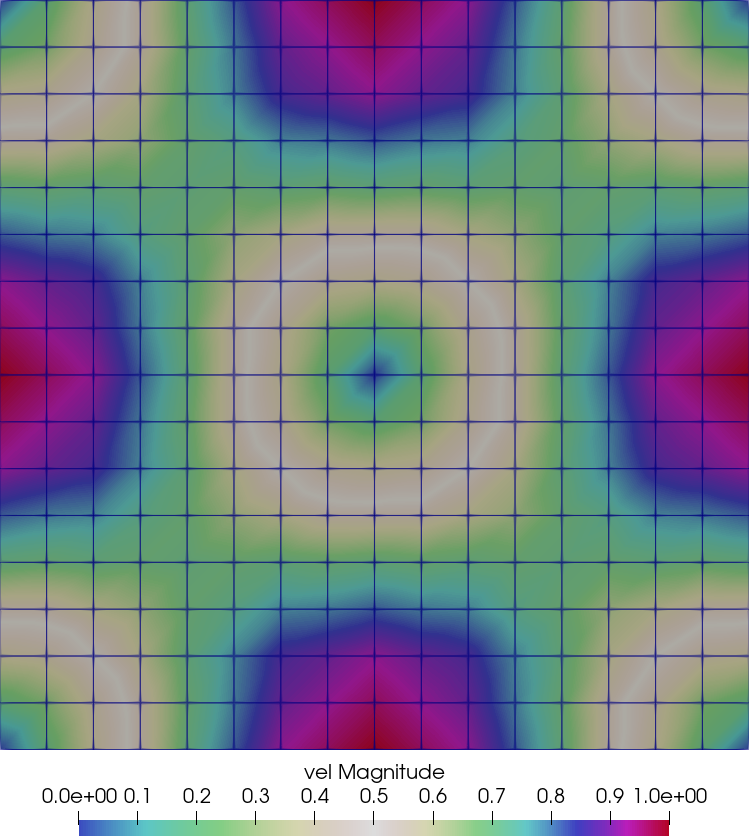}
		\end{subfigure}\hfill
		\begin{subfigure}[b]{0.19\textwidth}
			\includegraphics[width=\linewidth]{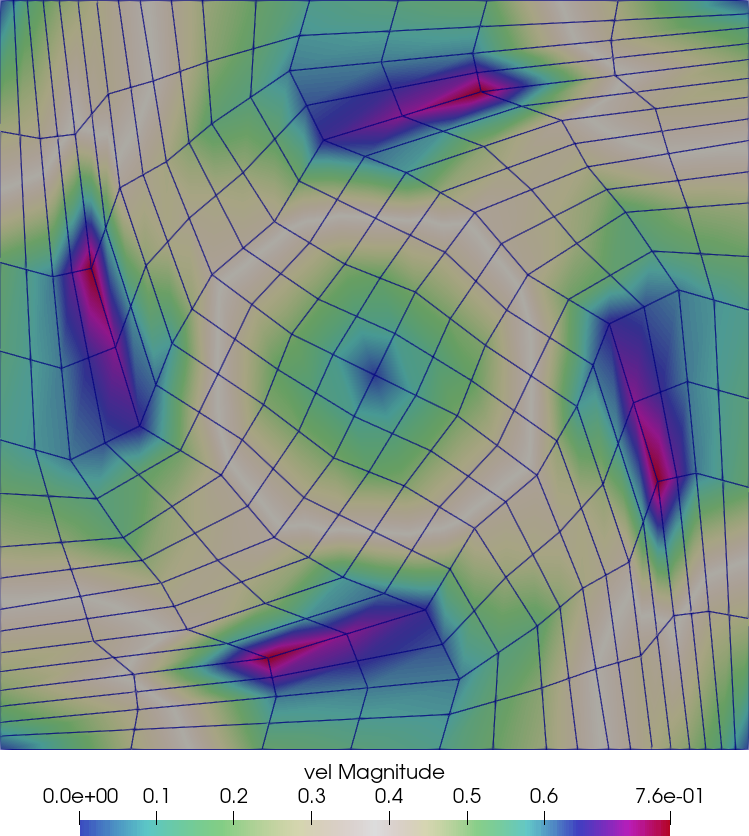}
		\end{subfigure}\hfill
		\begin{subfigure}[b]{0.19\textwidth}
			\includegraphics[width=\linewidth]{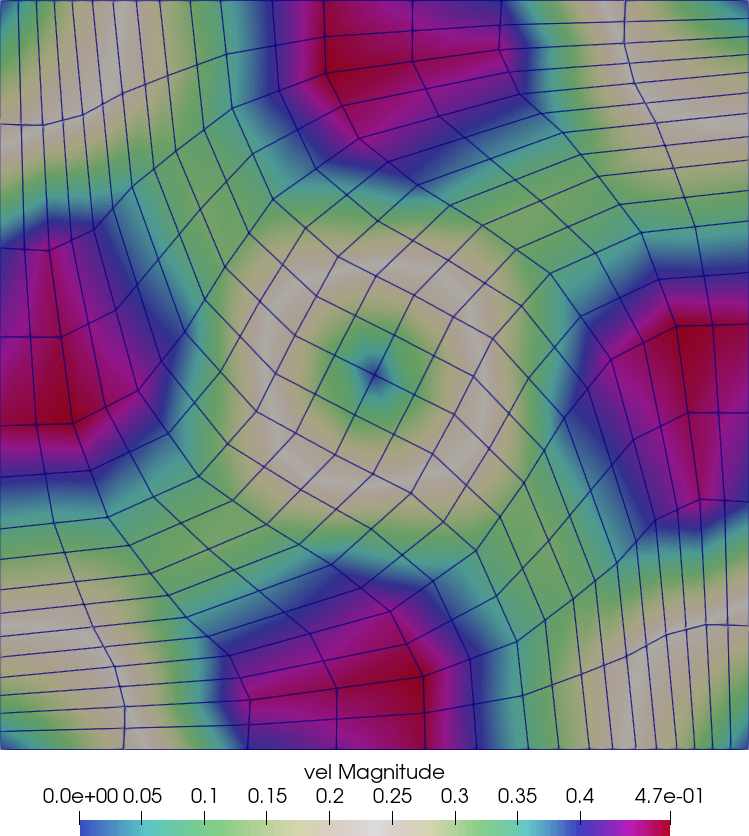}
		\end{subfigure}\hfill
		\begin{subfigure}[b]{0.19\textwidth}
			\includegraphics[width=\linewidth]{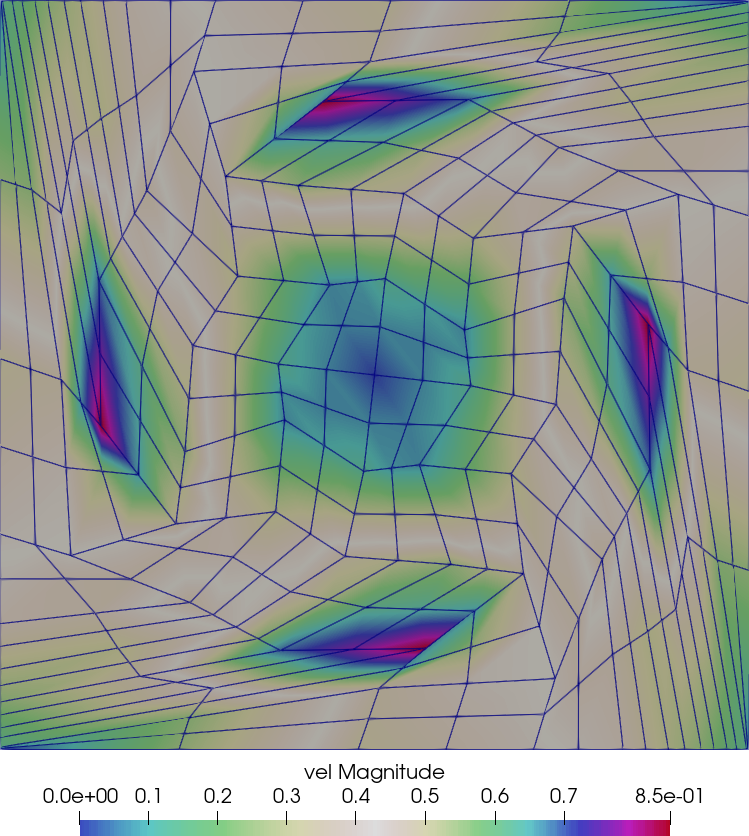}
		\end{subfigure}\hfill
		\begin{subfigure}[b]{0.19\textwidth}
			\includegraphics[width=\linewidth]{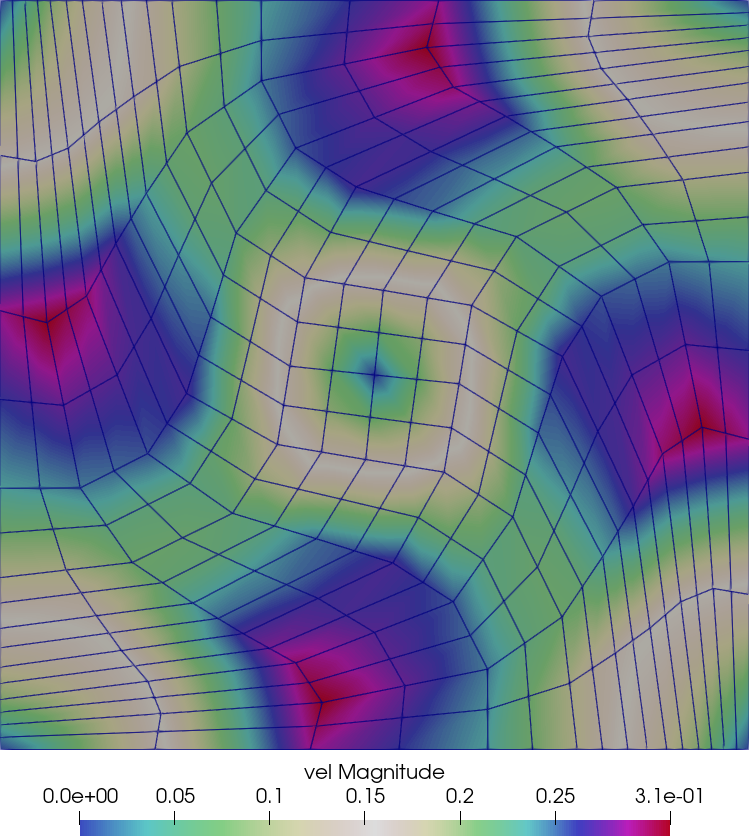}
		\end{subfigure}
		
		\vspace{6pt}

		\begin{subfigure}[b]{0.19\textwidth}
			\includegraphics[width=\linewidth]{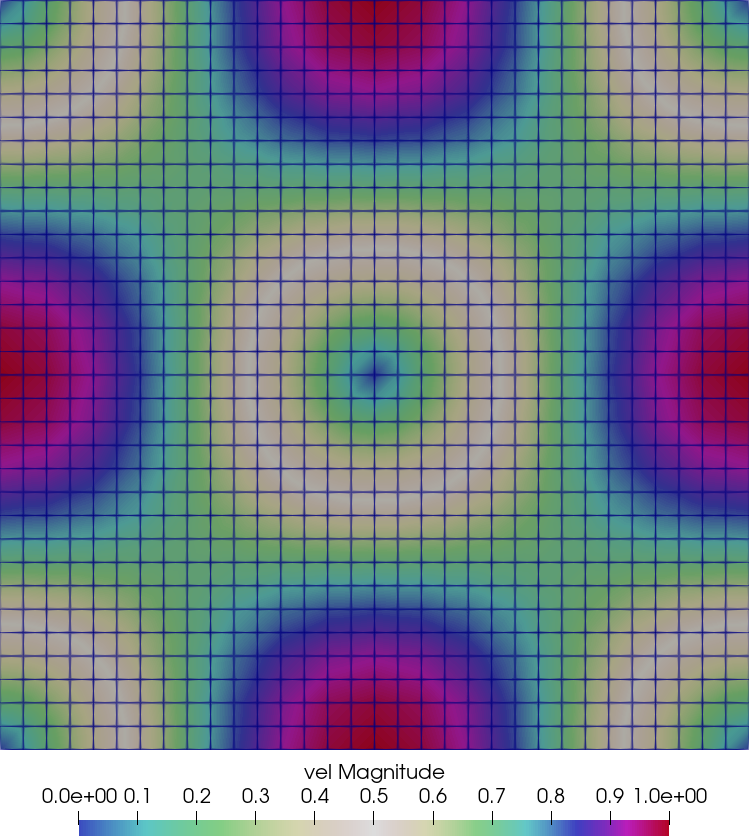}
		\end{subfigure}\hfill
		\begin{subfigure}[b]{0.19\textwidth}
			\includegraphics[width=\linewidth]{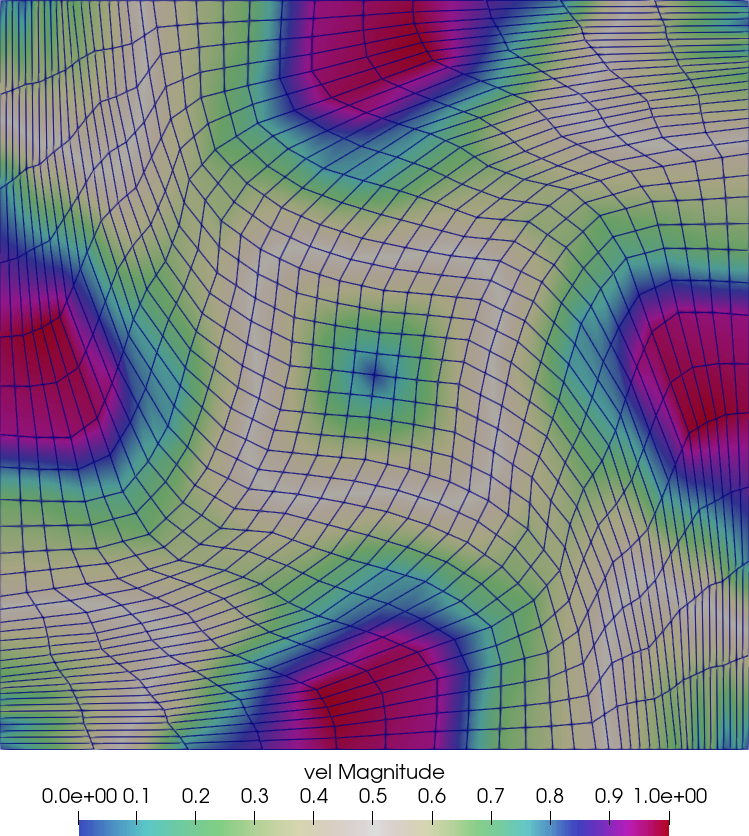}
		\end{subfigure}\hfill
		\begin{subfigure}[b]{0.19\textwidth}
			\includegraphics[width=\linewidth]{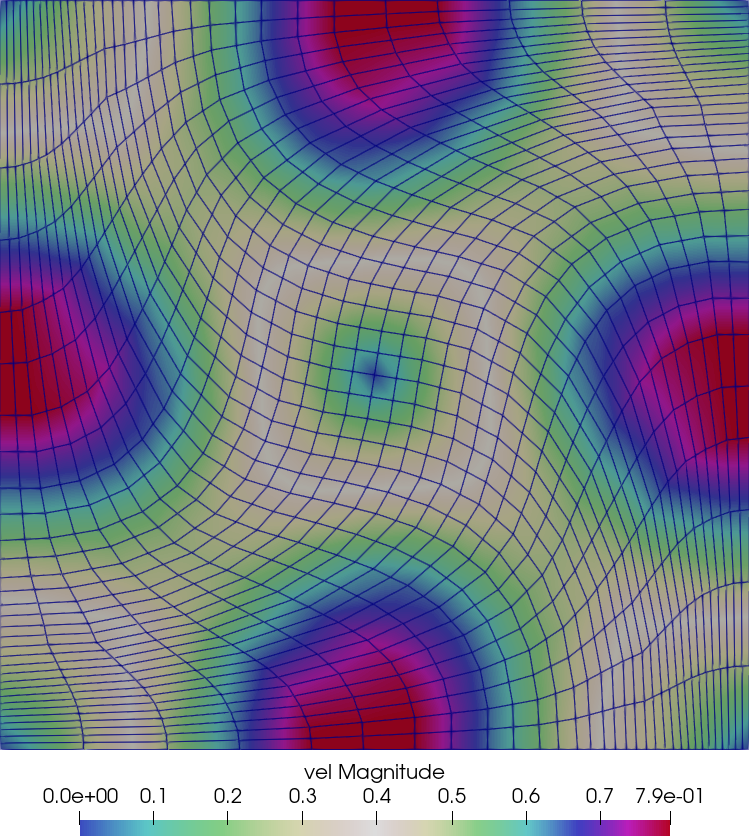}
		\end{subfigure}\hfill
		\begin{subfigure}[b]{0.19\textwidth}
			\includegraphics[width=\linewidth]{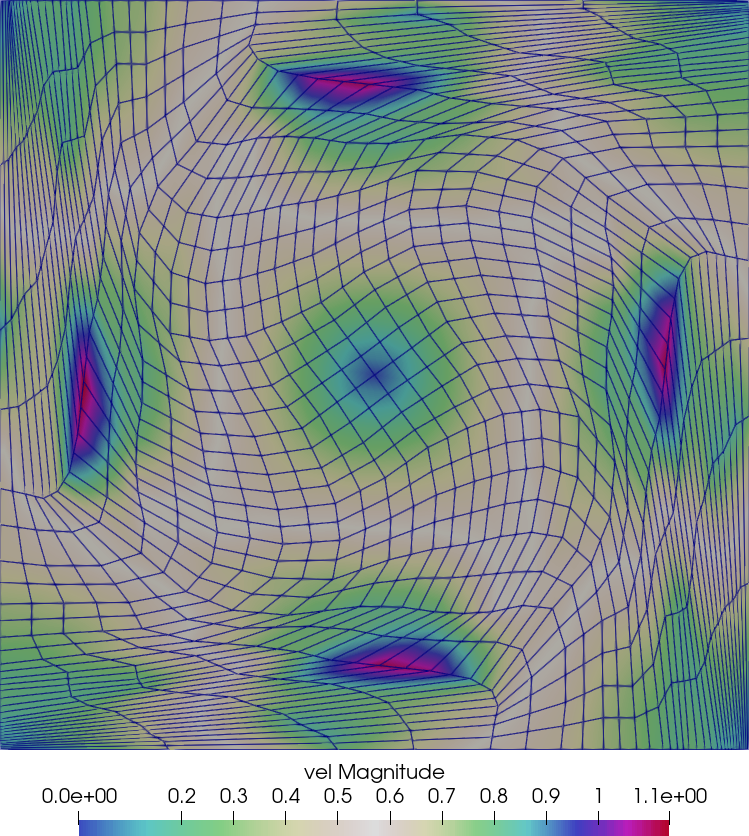}
		\end{subfigure}\hfill
		\begin{subfigure}[b]{0.19\textwidth}
			\includegraphics[width=\linewidth]{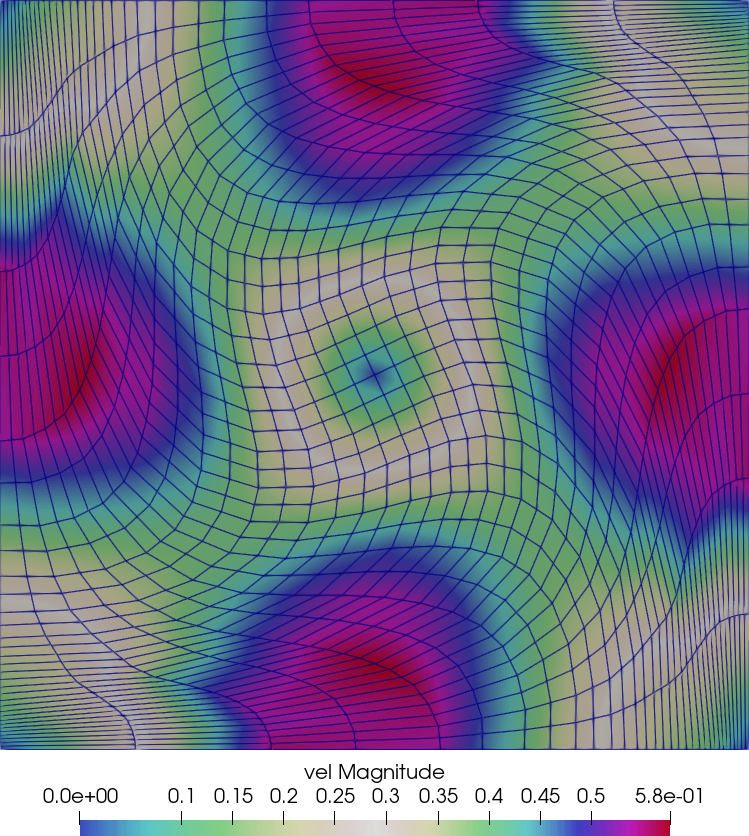}
		\end{subfigure}
		
		\vspace{6pt}
		
		\begin{subfigure}[b]{0.19\textwidth}
			\includegraphics[width=\linewidth]{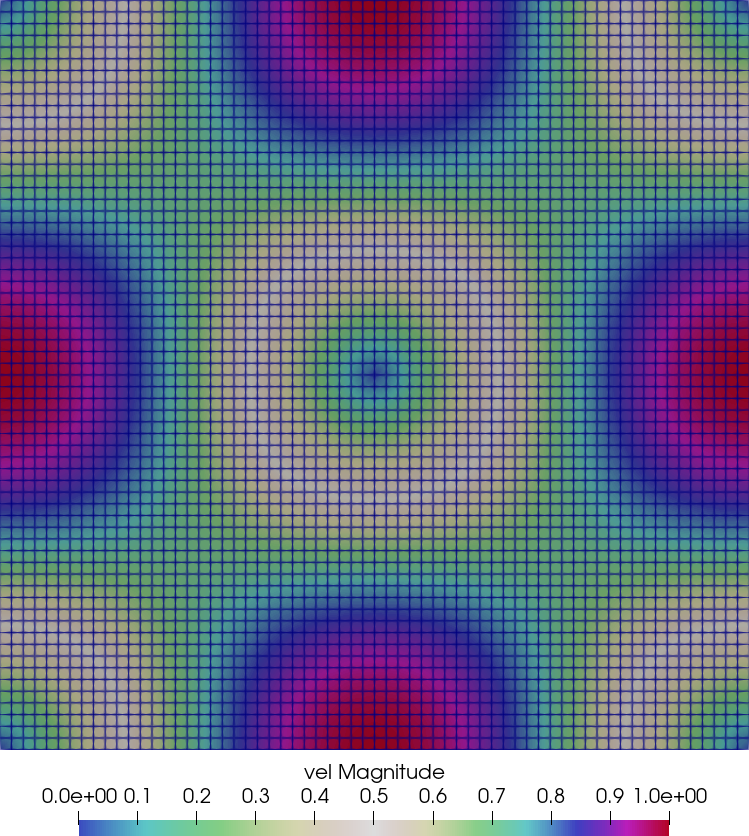}
		\end{subfigure}\hfill
		\begin{subfigure}[b]{0.19\textwidth}
			\includegraphics[width=\linewidth]{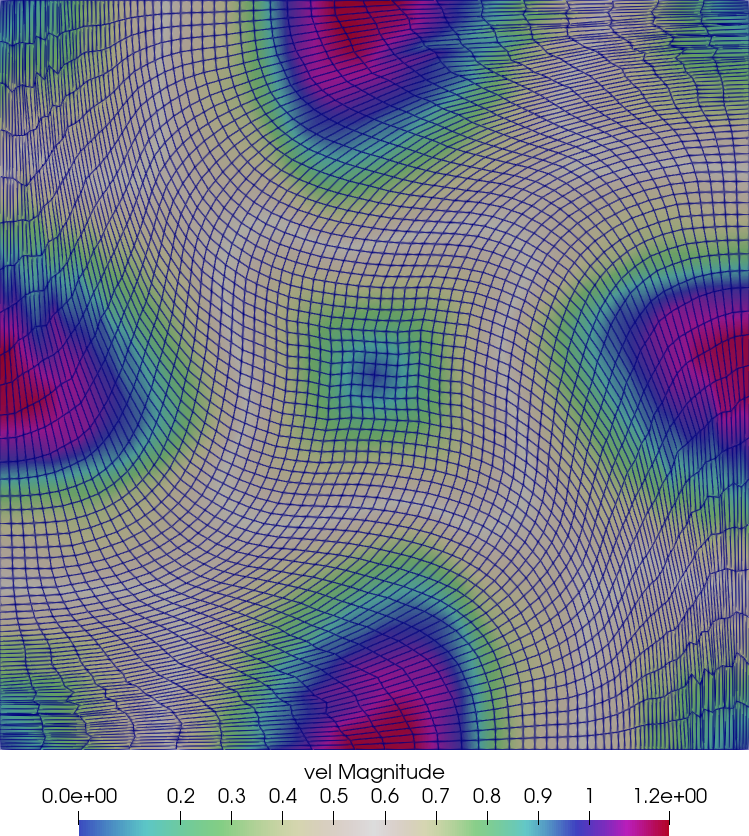}
		\end{subfigure}\hfill
		\begin{subfigure}[b]{0.19\textwidth}
			\includegraphics[width=\linewidth]{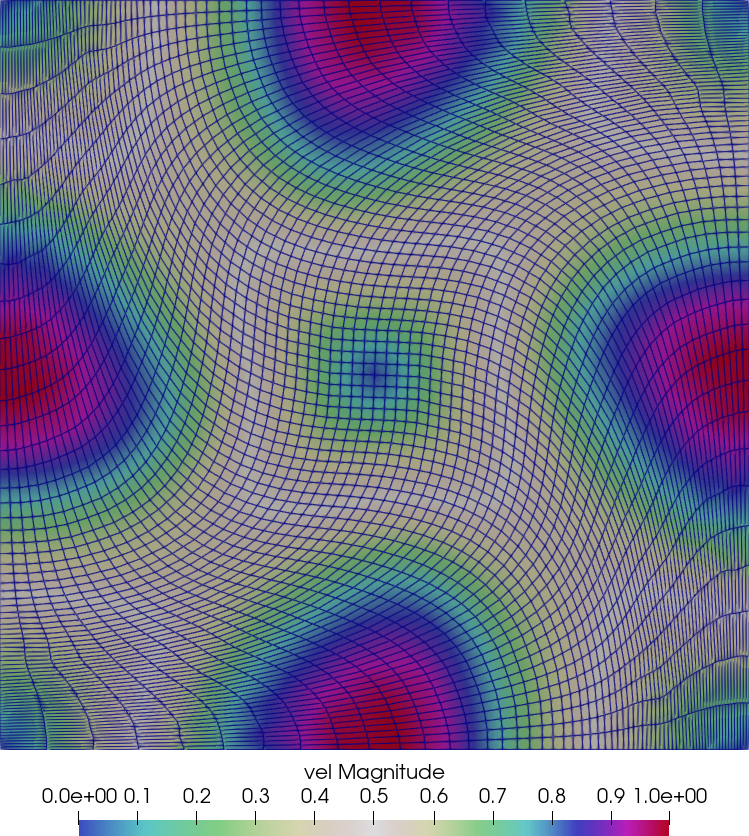}
		\end{subfigure}\hfill
		\begin{subfigure}[b]{0.19\textwidth}
			\includegraphics[width=\linewidth]{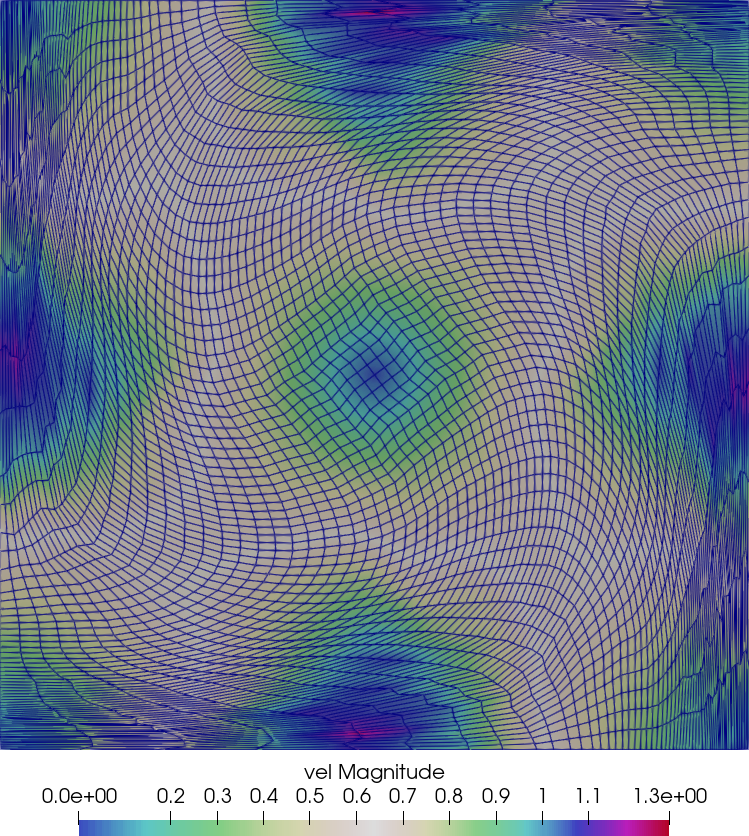}
		\end{subfigure}\hfill
		\begin{subfigure}[b]{0.19\textwidth}
			\includegraphics[width=\linewidth]{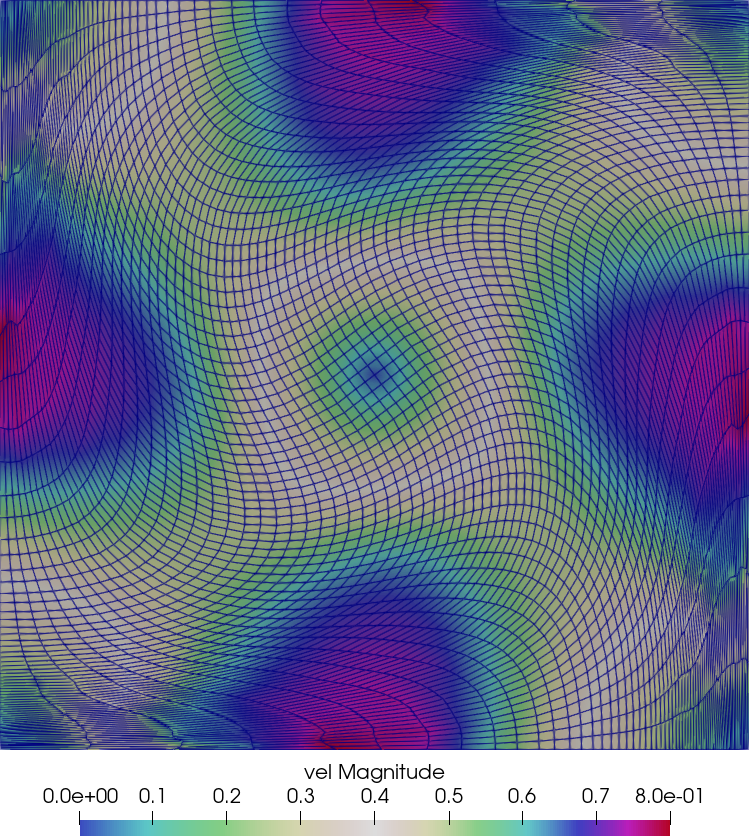}
		\end{subfigure}
		
		\vspace{6pt}
		
		\begin{subfigure}[b]{0.19\textwidth}
			\includegraphics[width=\linewidth]{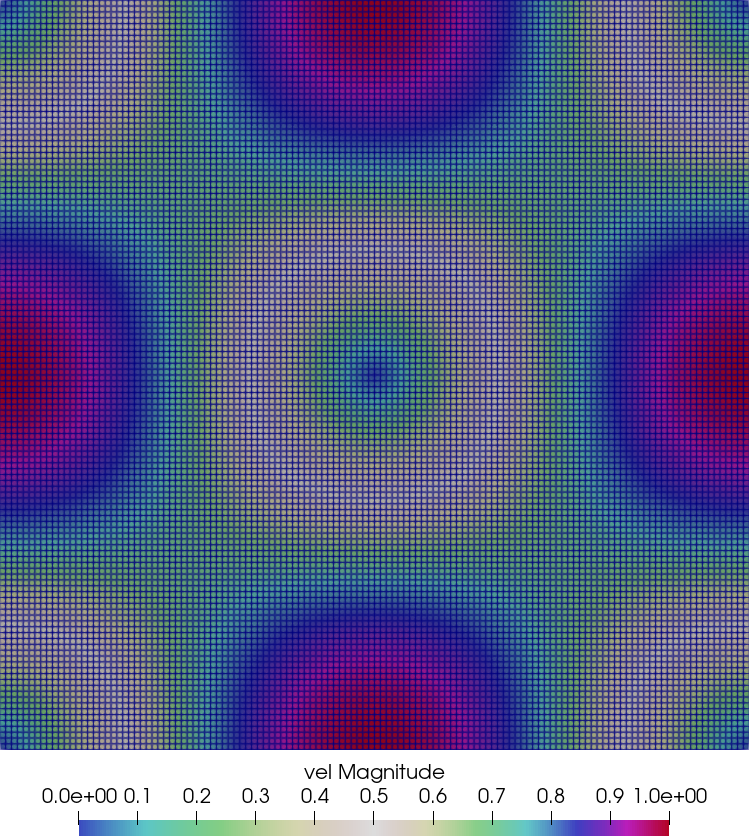}
			\caption{$t=0$}
		\end{subfigure}\hfill
		\begin{subfigure}[b]{0.19\textwidth}
			\includegraphics[width=\linewidth]{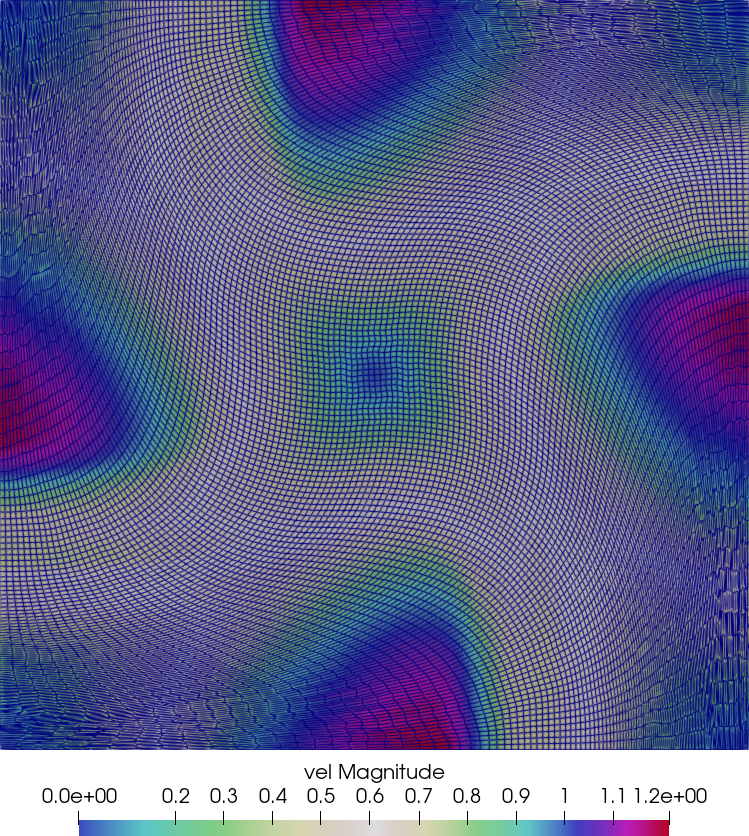}
			\caption{$t=0.5$, $\beta_c=1.0$}
		\end{subfigure}\hfill
		\begin{subfigure}[b]{0.19\textwidth}
			\includegraphics[width=\linewidth]{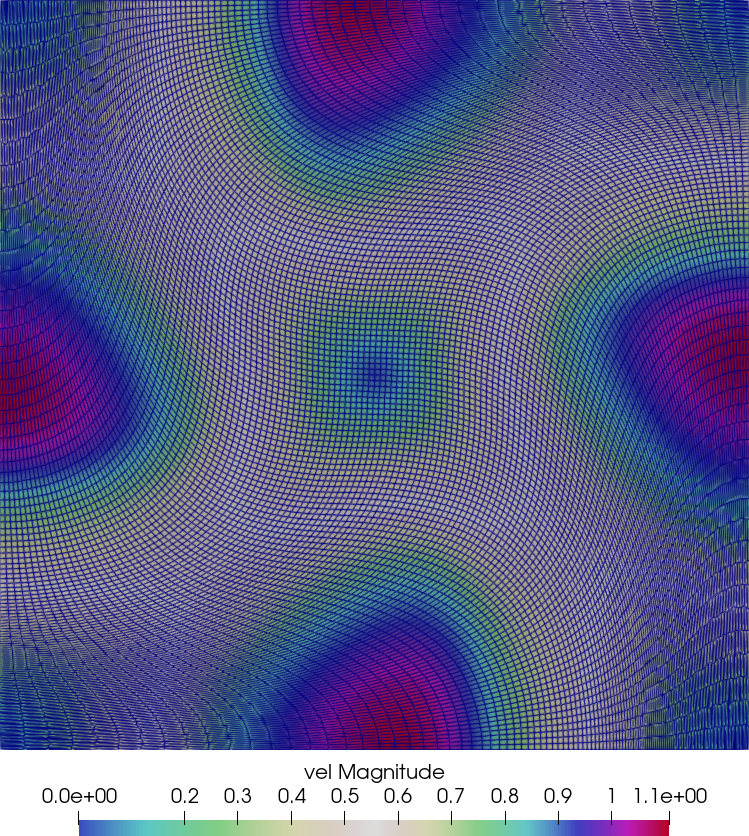}
			\caption{$t=0.5$, $\beta_c=0.8$}
		\end{subfigure}\hfill
		\begin{subfigure}[b]{0.19\textwidth}
			\includegraphics[width=\linewidth]{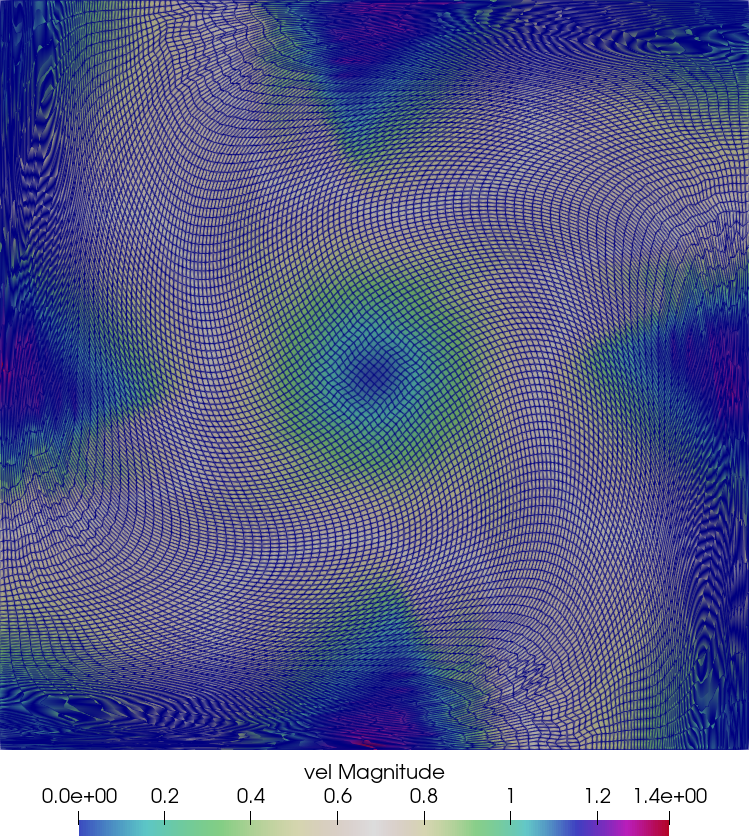}
			\caption{$t=0.75$, $\beta_c=1.0$}
		\end{subfigure}\hfill
		\begin{subfigure}[b]{0.19\textwidth}
			\includegraphics[width=\linewidth]{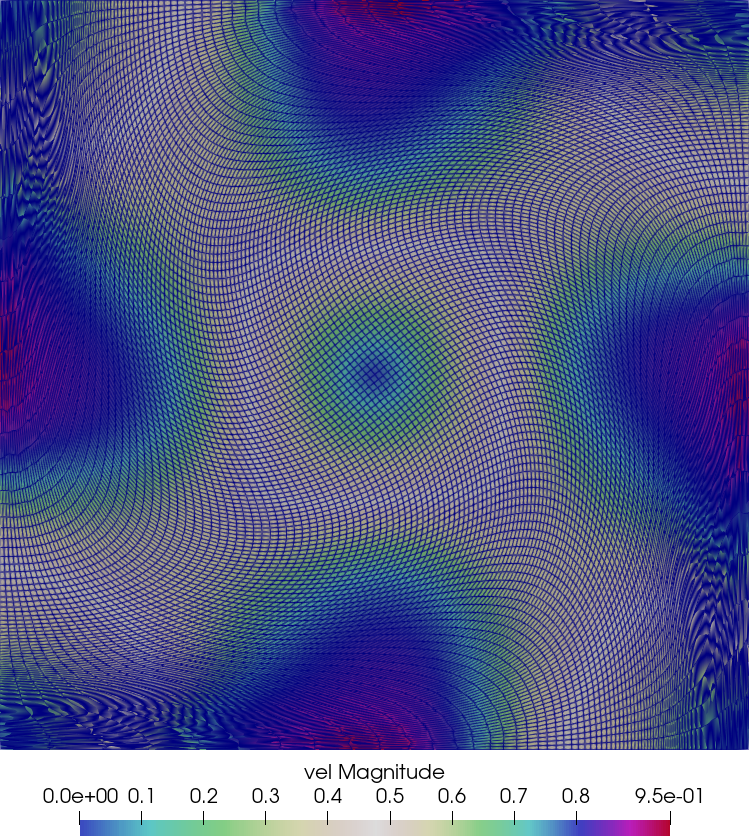}
			\caption{$t=0.75$, $\beta_c=0.8$}
		\end{subfigure}
		
		\caption{Taylor-Green vortex problem, velocity magnitude fields at different times $t$ for various mesh resolutions. Note that the solver outputs the subcells associated with each element.  See Fig.~\ref{fig:geom} for details. }
		\label{fig:tg}
	\end{figure}

	
	We conduct the simulation on uniformly refined hexahedral meshes consisting of a single-element-thick slab, with resolutions of $8 \times 8$, $16 \times 16$, $32 \times 32$, and $64 \times 64$. The velocity magnitude fields and corresponding meshes at times $t=0.0$, $t=0.5$ and $t=0.75$ are shown in Fig. \ref{fig:tg}. We observe mesh entangling when using $\beta_c=1.0$ which can be mitigated by reducing the parameter to $\beta_c=0.8$. Although the solution becomes severely distorted at later times, the simulation remains stable, and we are able to advance the computation beyond $t=2.0$ without any mesh remapping or code crashes, even when the mesh is heavily deformed.
	
	\subsection{3D Sedov blast wave problem}
	The Sedov problem \cite{pederson2016,sedov1959} involves an outward-traveling blast wave in a gamma-law gas generated by an energy source. We set $\gamma=7/5$ and test our method for the three-dimensional case. An exact solution, derived from self-similarity arguments, is available; see, e.g., \cite{kamm2007}. The initial conditions read
	\begin{align*}
		\varrho^0=1, \quad \mathbf{u}^0=0, \quad p^0=10^{-6}.
	\end{align*}
	The energy source is positioned at the origin, and the pressure in the element containing the origin is calculated as $p_{or}=(\gamma-1)\varrho_{or}\frac{e_0}{v_{or}}$, where $e_0$ denotes the total amount of released internal energy and $v_{or}$ is the volume of the element containing the origin. By setting $e_0=0.493390$, the solution features a diverging shock wave of infinite strength, with the front located at a radius of $r=1$ at $t=1$, and a peak density reaching $6$.  The computational domain is $\Omega=(0,1.2)^3$.
	
\begin{figure}
		\centering
		\begin{subfigure}[b]{0.32\textwidth}
			\hspace{1cm}
			\includegraphics[width=\linewidth]{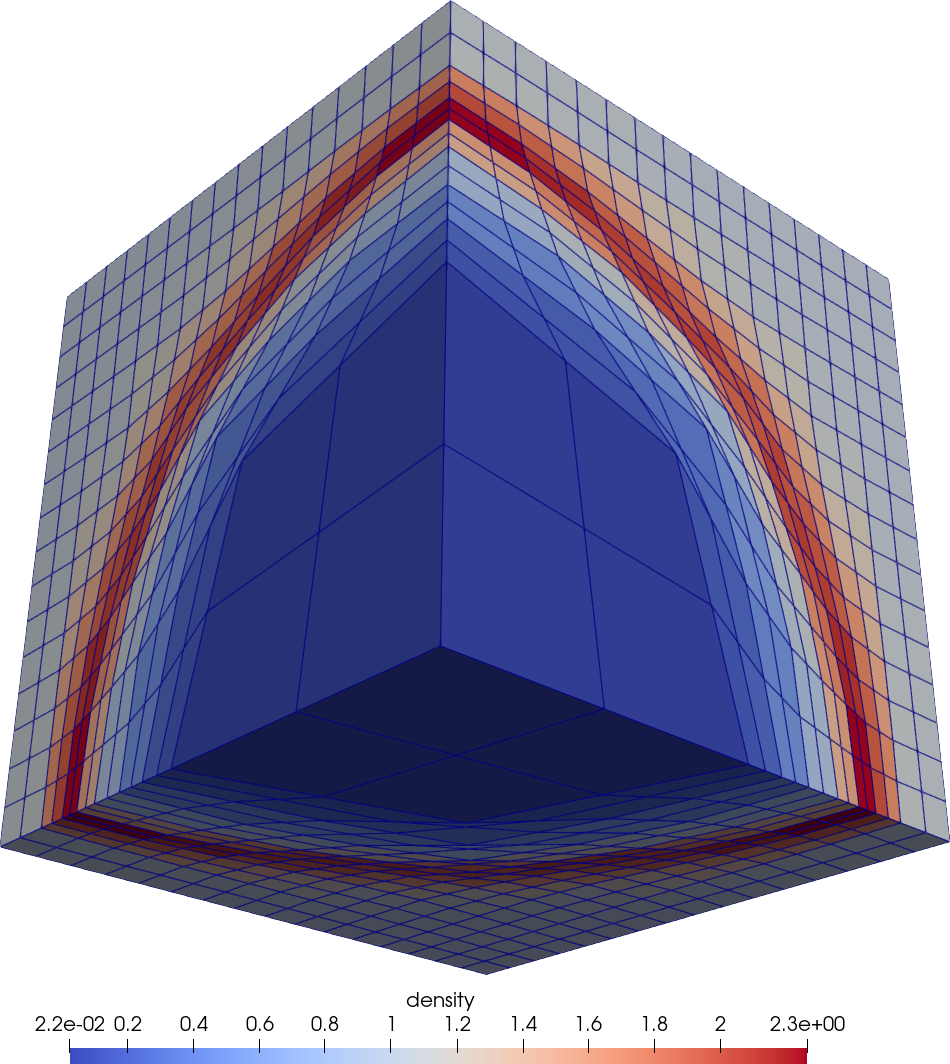}
		\end{subfigure}\hfill
		\begin{subfigure}[b]{0.45\textwidth}
			\includegraphics[width=\linewidth]{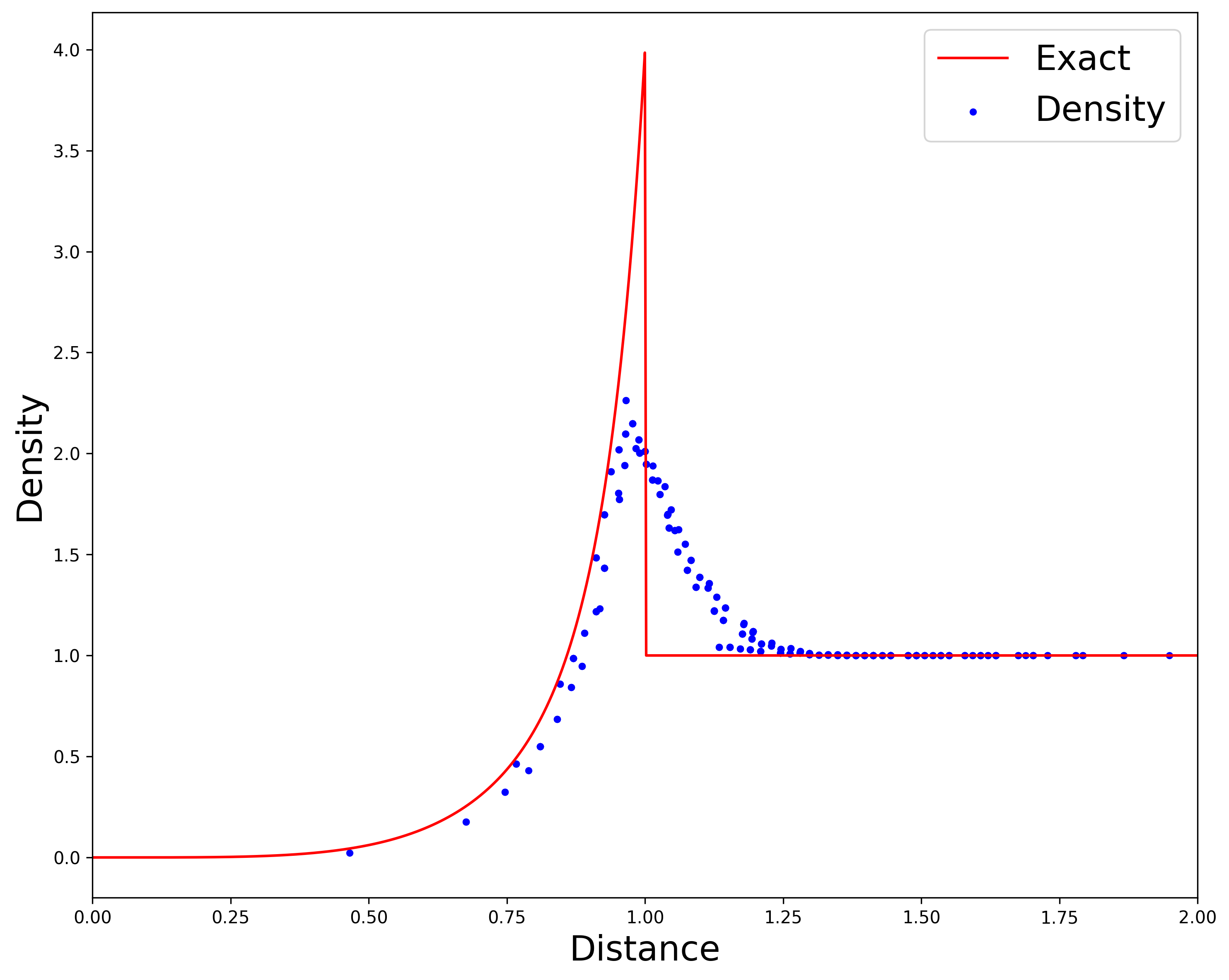}
		\end{subfigure}
		
		\vspace{6pt}
		
		\begin{subfigure}[b]{0.32\textwidth}
			\hspace{1cm}
			\includegraphics[width=\linewidth]{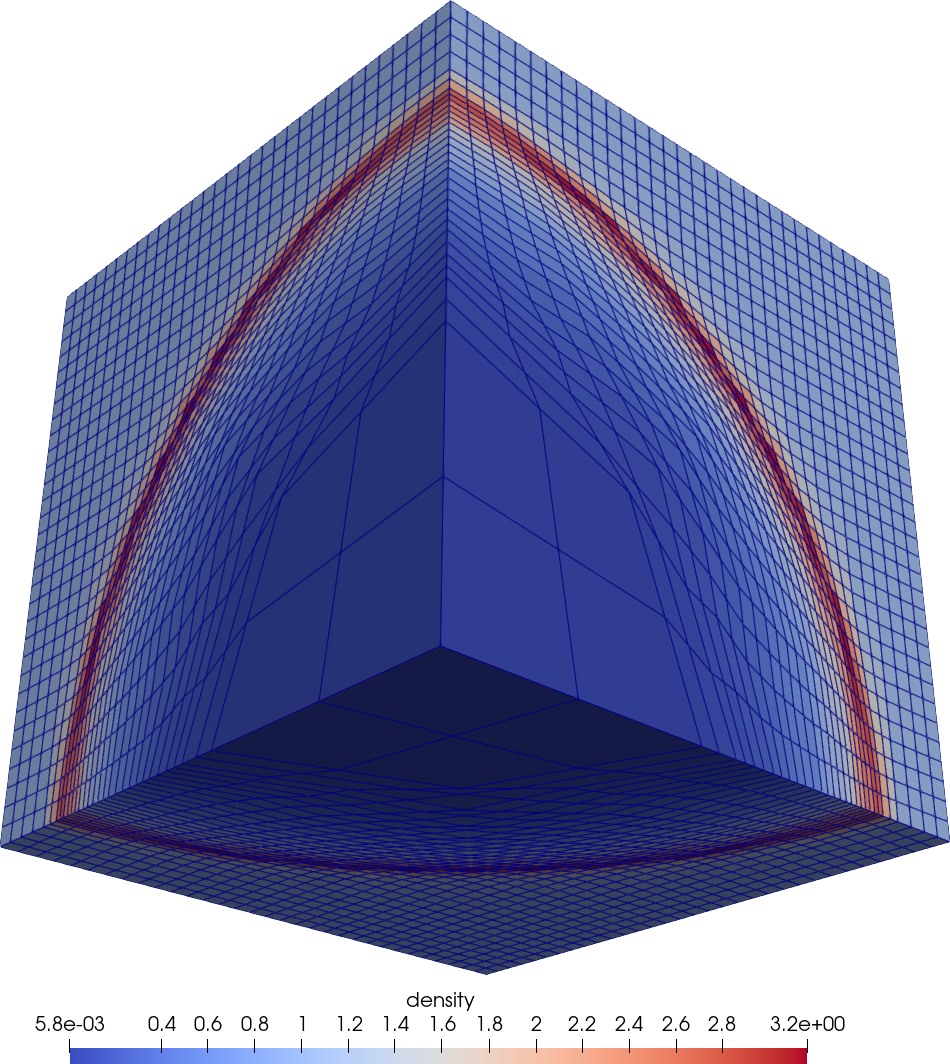}
		\end{subfigure}\hfill
		\begin{subfigure}[b]{0.45\textwidth}
			\includegraphics[width=\linewidth]{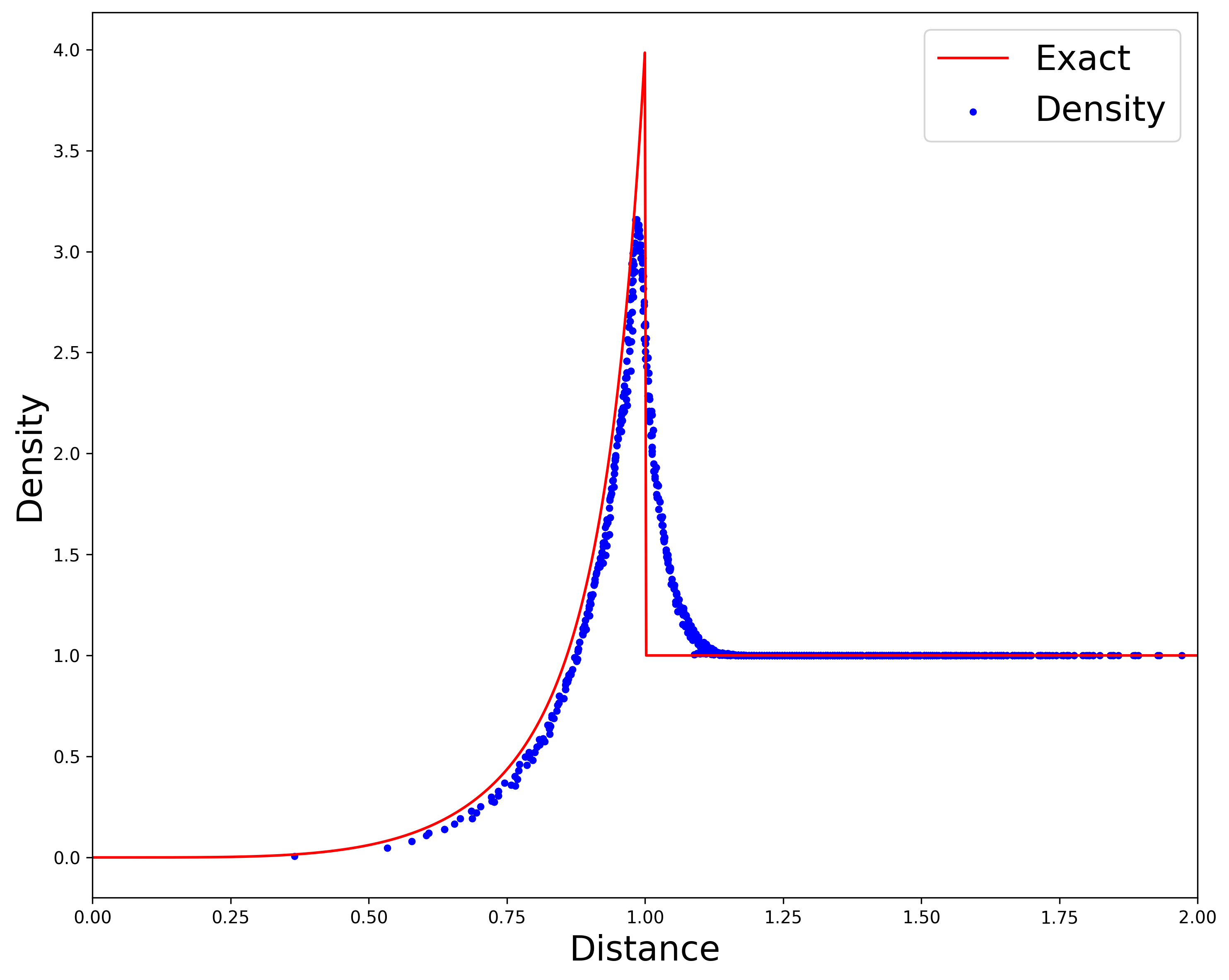}
		\end{subfigure}
		
		\vspace{6pt}
		
		\begin{subfigure}[b]{0.32\textwidth}
			\hspace{1cm}
			\includegraphics[width=\linewidth]{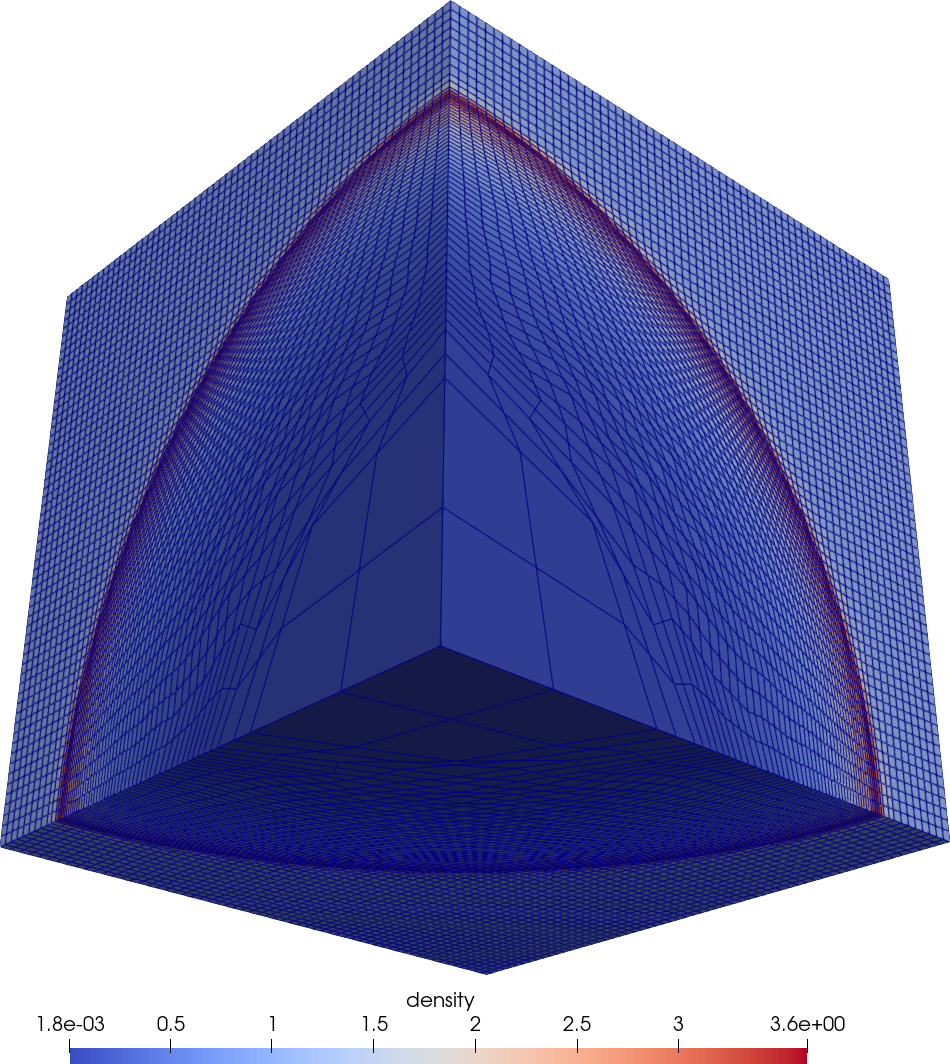}
		\end{subfigure}\hfill
		\begin{subfigure}[b]{0.45\textwidth}
			\includegraphics[width=\linewidth]{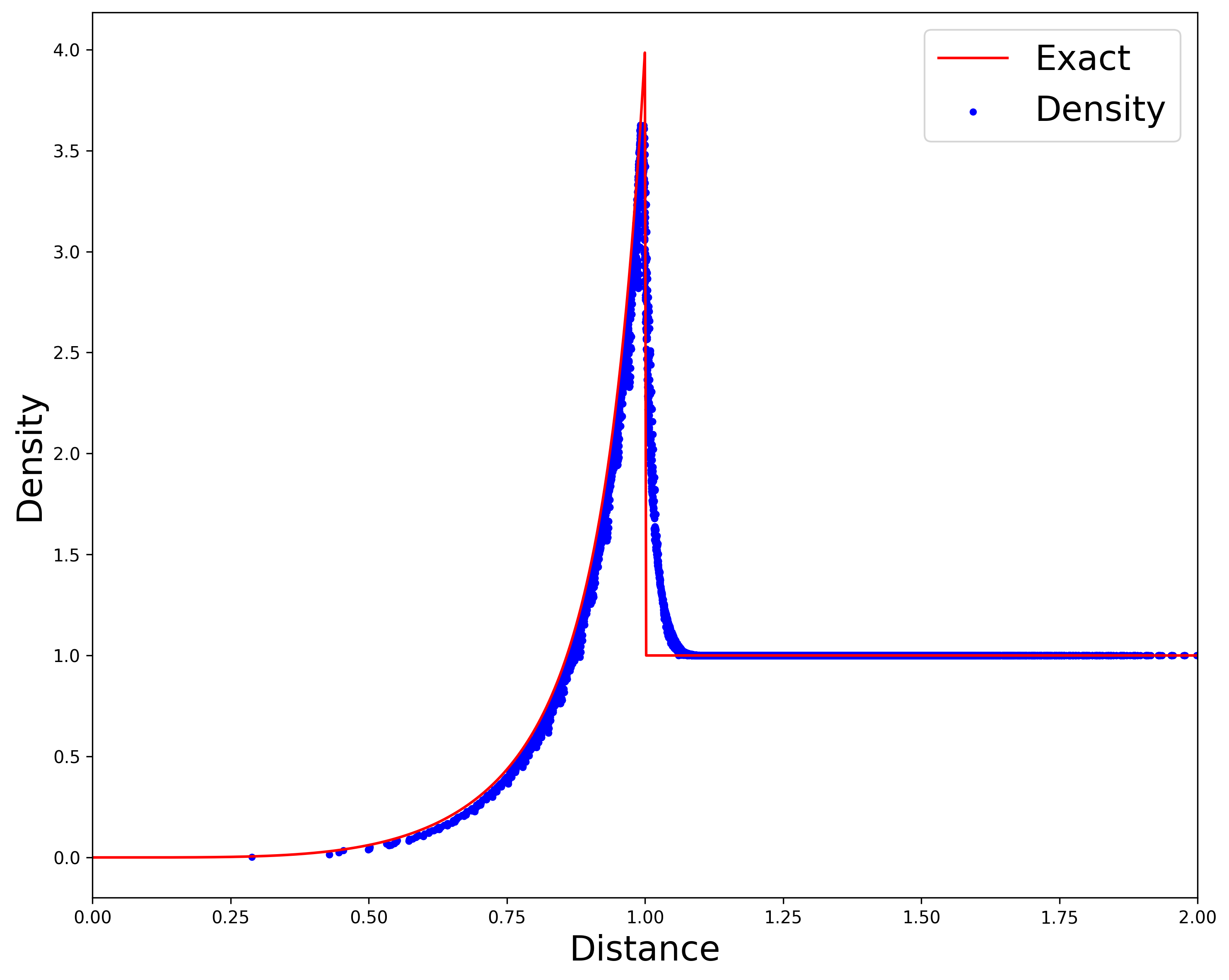}
		\end{subfigure}
		
		\caption{Three-dimensional Sedov blast wave problem, density profiles $\varrho$ at $t=1.0$ for various mesh resolutions. Note that the solver outputs the subcells associated with each element.  See Fig.~\ref{fig:geom} for details.}
		\label{fig:sed}
	\end{figure}

	We split the computational domain into hexahedral meshes with $8$, $16$, and $30$ elements in each spatial direction. Fig. \ref{fig:sed} shows the density profiles and corresponding scatter plots of element-averaged density versus radius. For all mesh resolutions, the shock front is accurately located and the solution maintains cylindrical symmetry. As the mesh is refined, the density peak nearly reaches $6$.
	
	\subsection{3D Noh problem}
	The Noh problem \cite{noh1987} is a standard test case used to validate the accuracy and symmetry preservation of Lagrangian methods. In this scenario, a cold gas with $\gamma=5/3$ is initially given an inward radial velocity, generating a diverging cylindrical shock wave that propagates at a speed of $s=1/3$. The density plateau behind the shock wave reaches a value of 16. The initial flow field is given by
	\begin{align*}
		\varrho^0=1, \quad \mathbf{u}^0=(-x/\sqrt{x^2+y^2},-y/\sqrt{x^2+y^2})^\top, \quad p^0=10^{-6}.
	\end{align*}
	The exact solution is an infinite-strength symmetric shock; see, e.g., \cite{dobrev2018,liska2003}. The initial computational domain is $\Omega^0=(0,1)\times(0,1)$, with wall boundary conditions applied along the lines $x=0$ and $y=0$. We run the simulation until $t=0.6$. As noted by Liska and Wendroff \cite{liska2003}, many high-resolution schemes tend to produce highly oscillatory results or fail altogether.
	
\begin{figure}
		\centering
		\begin{subfigure}[b]{0.48\textwidth}
			\includegraphics[width=\linewidth]{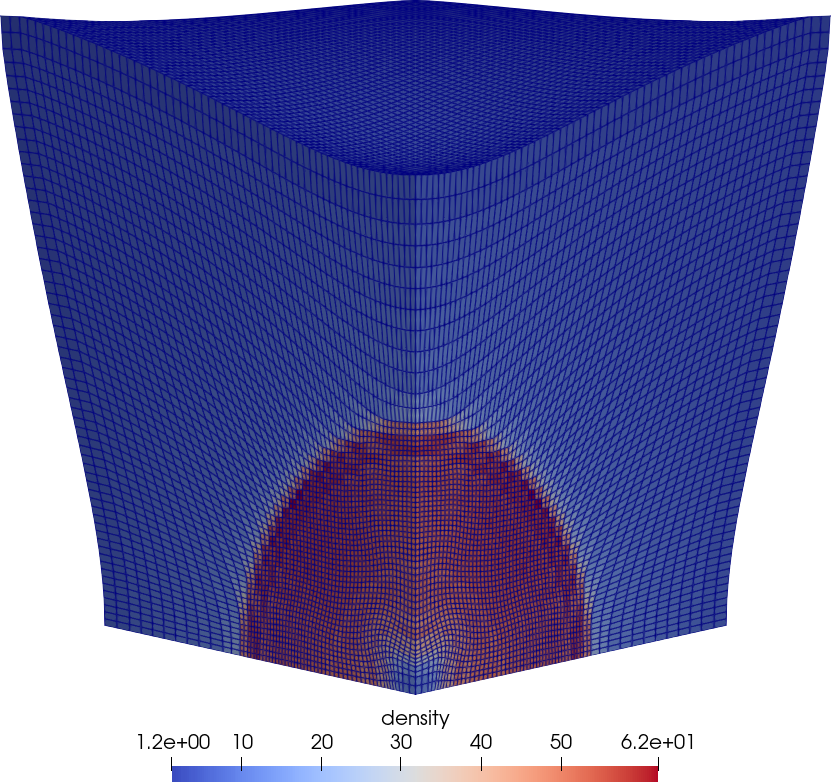}
			\caption{$\beta_c=0.8$}
		\end{subfigure}\hfill
		\begin{subfigure}[b]{0.48\textwidth}
			\includegraphics[width=\linewidth]{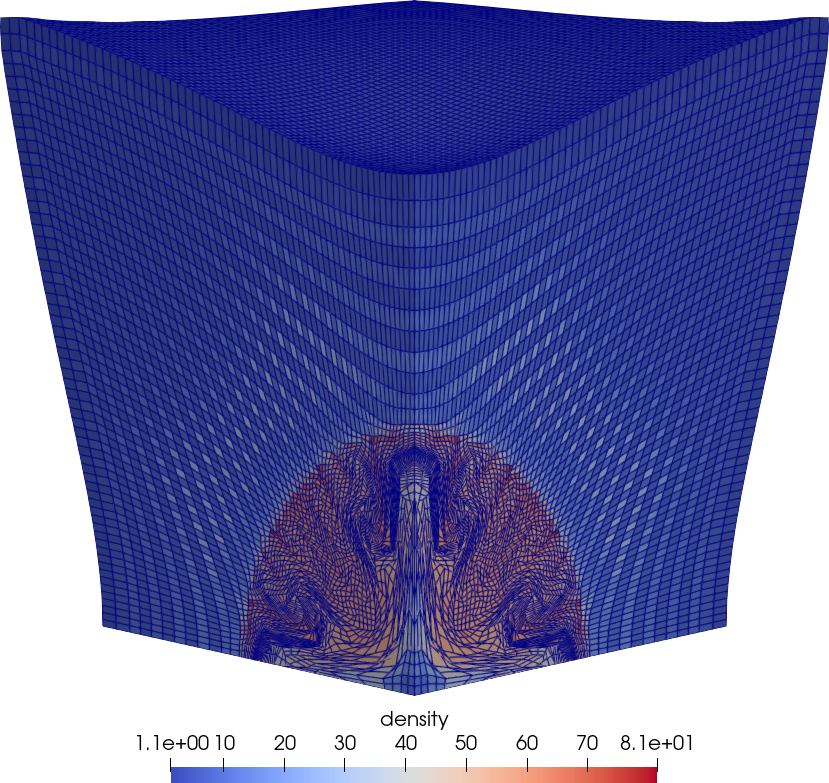}
			\caption{$\beta_c=1.0$}
		\end{subfigure}
		\caption{Noh problem, density profiles $\varrho$ obtained at $t=0.6$ on a $30\times 30\times 30$ mesh. Note that the solver outputs the subcells associated with each element.  See Fig.~\ref{fig:geom} for details.}
		\label{fig:noh}
	\end{figure}

	We simulate the Noh problem on a uniform hexahedral mesh of size $30\times 30 \times 30$, making this a challenging test case since the mesh is not aligned with the flow. The density profiles obtained with our DGH method at time $t=0.6$ are shown in Fig. \ref{fig:noh}. Mesh entangling occurs for $\beta_c=1.0$, but can be significantly reduced by using $\beta_c=0.8$. Notably, the DGH method does not suffer from IDP violations even for highly distorted meshes. The shock is located on a circle of radius approximately $0.2$, consistent with results from the literature \cite{vilar2012a,vilar2014}. For $\beta_c=0.8$, the shock wave front is sharply resolved and closely matches the one-dimensional cylindrical solution, while the method successfully preserves the radial symmetry of the flow.
	
	\subsection{2D Triple point problem}
	In our final numerical experiment, we consider the triple point problem \cite{galera2010}, a three-state, two-material 2D Riemann problem known for generating vorticity. This problem is commonly used to test the robustness of Lagrangian methods in handling complex phenomena, such as significant vorticity, large shear, and interacting shocks \cite{loubere2005}. The initial setup consists of three regions: a high-pressure region on the left and two low-pressure regions on the right. The high-pressure region drives a shock wave through the adjacent low-pressure regions, leading to the formation of a vortex at the triple point. The initial configuration on the computational domain $\Omega=(0,7)\times(0,3)$ is given by
	\begin{align*}
		(\varrho^0,\mathbf{u}^0,p^0)^\top=
		\begin{cases}
			(1,0,1)^\top \; &\text{if} \; y\leq 1,\\
			(0.1,0,0.1)^\top \; &\text{if} \; x \geq 1, y \geq 1.5, \\
			(1,0,0.1)^\top \; &\text{if} \; x \geq 1, y < 1.5.
		\end{cases}
	\end{align*}
	We set $\gamma=1.5$ in the initial left and upper regions, and $\gamma=1.4$ in the lower region.
	
\begin{figure}
		\centering
		\begin{subfigure}[b]{0.48\textwidth}
			\includegraphics[width=\linewidth]{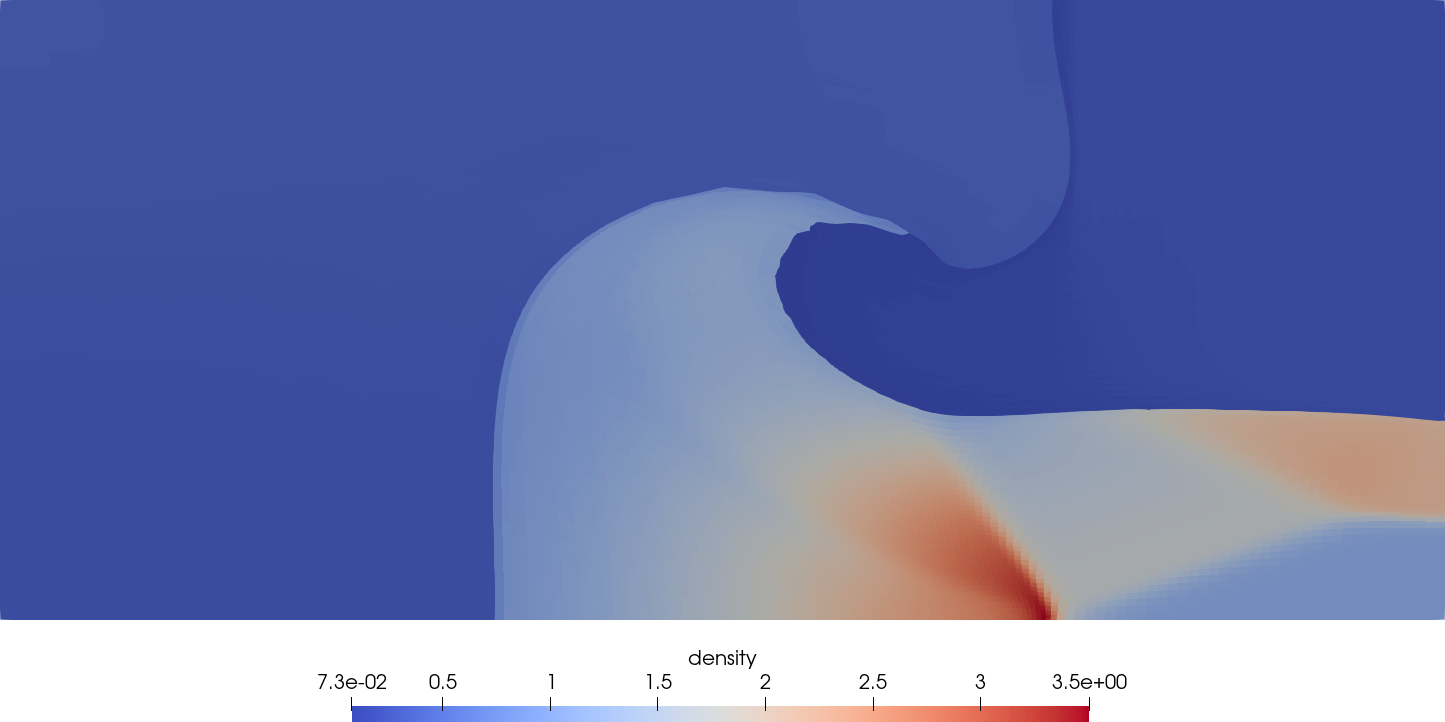}
		\end{subfigure}\hfill
		\begin{subfigure}[b]{0.48\textwidth}
			\includegraphics[width=\linewidth]{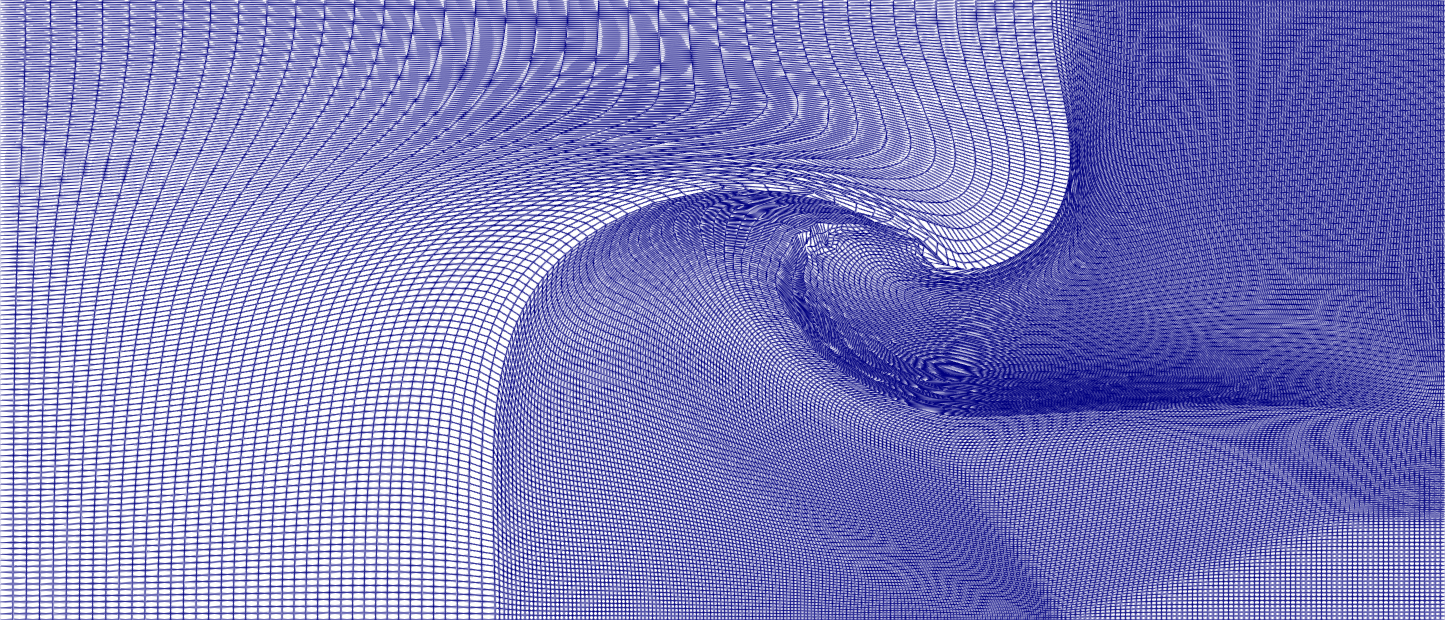}
			\vspace{0.1cm}
		\end{subfigure}
		\caption{Triple point problem, density profile $\varrho$ and mesh at time $t=4.0$ using $E_h=140\times 80$ elements. Note that the solver outputs the subcells associated with each element.  See fig.\ref{fig:geom} for details.}
		\label{fig:tpp}
	\end{figure}

	We run the simulation using a single-element-thick slab with $E_h=140\times 80$ uniform elements and $\beta_c=0.8$. The density profile and computational mesh at $t=4.0$ are presented in Fig. \ref{fig:tpp}. The solution demonstrates the ability of our DGH method to accurately capture sharp interfaces between the interacting shock waves. The triple point region is accurately resolved, clearly displaying the physical vorticity.
	
	\section{Conclusions}
	\label{sec:concl}
	We discussed a new way to limit Lagrangian DG discretizations of the Euler equations. The proposed scheme is second-order accurate while remaining fully conservative and GCL consistent. Following FCT design principles, we combine a property-preserving low-order scheme with a high-order target discretization. By utilizing a nodal version of Zalesak's limiting strategy, we achieve optimal convergence rates and effectively suppress spurious oscillations globally. To address local oscillations, we employ the clip-and-scale slope limiter, which yields a locally bound-preserving conservative approximation. Our analysis highlights the importance of using second-order time integrators to maintain GCL consistency. Additionally, it is worthwile to explore FCT-type limiters for systems, where correction factors are computed individually for each field, as opposed to the scalar limiting approach we used. Extending our approach to higher-order finite elements and curvilinear elements offers another potential avenue for further investigation.

    \section*{Acknowledgments}
    Nathaniel Morgan and Jacob Moore gratefully acknowledge the funding from the Laboratory Directed Research and Development (LDRD) program at Los Alamos National Laboratory (LANL) under project number 20230028DR. The Los Alamos unlimited release number is LA-UR-25-30128. LANL is operated by Triad National Security, LLC for the U.S. Department of Energy’s National Nuclear Security Administration under contract number 89233218CNA000001.
    
	\appendix
	\section*{Appendix. GCL consistency for triangles and quadrilaterals.}
	\label{sec:appendix}
	We consider a triangle with vertices $\mathbf{x}_p=(x_p,y_p)$, $p\in \mathcal{P}(c)$, numbered counterclockwise. Each node moves with a constant velocity $\mathbf{u}_p=(u_p^x,u_p^y)^\top$ over the time interval $[t^n,t^{n+1}]$. The nodal coordinates are defined as functions of time:
	\begin{align}
	x_p(t)=x_p^n+u_p^x(t-t^n), \qquad y_p(t)=y_p^n+u_p^y(t-t^n), \qquad t\in[t^n,t^{n+1}].
	\label{eq:coordupd}
	\end{align}
	Inserting $t=t^{n+1}$ into \eqref{eq:coordupd} yields
	\begin{align*}
	x_p^{n+1}=x_p^n+\Delta t u_p^x, \qquad y_p^{n+1}=y_p^n+\Delta t u_p^y. 
	\end{align*}
	
	The algebraic volume $\tilde{v}_c$ at time $t^{n+1}$ can be computed using the discretized GCL
	\begin{align}
	\tilde{v}_c^{n+1}=\tilde{v}_c^n+\Delta t\sum_{p\in \mathcal{P}(c)}\sum_{s_f\in \mathcal{SF}(c,p)}a_{cps_f}\mathbf{u}_p\cdot \mathbf{n}_{cps_f}=\tilde{v}_c^n+\sum_{p\in \mathcal{P}(c)}\sum_{s_f\in \mathcal{SF}(c,p)}a_{cps_f}(\mathbf{x}_p^{n+1}-\mathbf{x}_p^n)\cdot \mathbf{n}_{cps_f}.
	\label{eq:algvol}
	\end{align} 
	The change in algebraic volume follows as
	\begin{align*}
	&\tilde{v}_c^{n+1}-\tilde{v}_c^n=\sum_{p\in \mathcal{P}(c)}\sum_{s_f\in \mathcal{SF}(c,p)}a_{cps_f}(\mathbf{x}_p^{n+1}-\mathbf{x}_p^n)\cdot \mathbf{n}_{cps_f}\\&=\frac{1}{2}[(x_1^{n+1}-x_1^n)(y_1^n-y_3^n)+(y_1^{n+1}-y_1^n)(x_3^n-x_1^n)
	+(x_1^{n+1}-x_1^n)(y_2^n-y_1^n)+(y_1^{n+1}-y_1^n)(x_1^n-x_2^n)\\
	&+(x_2^{n+1}-x_2^n)(y_2^n-y_1^n)+(y_2^{n+1}-y_2^n)(x_1^n-x_2^n)
	+(x_2^{n+1}-x_2^n)(y_3^n-y_2^n)+(y_2^{n+1}-y_2^n)(x_2^n-x_3^n)\\
	&+(x_3^{n+1}-x_3^n)(y_3^n-y_2^n)+(y_3^{n+1}-y_3^n)(x_2^n-x_3^n)
	+(x_3^{n+1}-x_3^n)(y_1^n-y_3^n)+(y_3^{n+1}-y_3^n)(x_3^n-x_1^n)]\\
	&=\frac{1}{2}[(x_1^{n+1}-x_1^n)(y_2^n-y_3^n)+(y_1^{n+1}-y_1^n)(x_3^n-x_2^n)
	+(x_2^{n+1}-x_2^n)(y_3^n-y_1^n)+(y_2^{n+1}-y_2^n)(x_1^n-x_3^n)\\
	&+(x_3^{n+1}-x_3^n)(y_1^n-y_2^n)+(y_3^{n+1}-y_3^n)(x_2^n-x_1^n)].
	\end{align*}
	
	The geometric volume $v_c$ of the triangle is
	\begin{align*}
	v_c=\frac{1}{2}\det \begin{pmatrix}
	1 & x_1 & y_1 \\ 1 & x_2 & y_2 \\ 1 & x_3 & y_3
	\end{pmatrix}
	=\frac{1}{2}(x_1y_2-y_1x_2+y_1x_3-x_1y_3+x_2y_3-x_3y_2).
	\end{align*}
	The change in geometric volume can be written as \cite{loubere2008}
	\begin{align*}
	v_c^{n+1}-v_c^n=\sum_{p=1}^{3}u_p^x\int_{t^n}^{t^{n+1}}\frac{\partial v_c}{\partial x_p}\, \mathrm{d}t+\sum_{p=1}^{3}u_p^y\int_{t^n}^{t^{n+1}}\frac{\partial v_c}{\partial y_p}\, \mathrm{d}t.
	\end{align*}
	We have
	\begin{alignat*}{3}
	\frac{\partial v_c}{\partial x_1}&=\frac{1}{2}(y_2-y_3), \qquad \frac{\partial v_c}{\partial x_2}&=\frac{1}{2}(y_3-y_1), \qquad \frac{\partial v_c}{\partial x_3}&=\frac{1}{2}(y_1-y_2), \\
	\frac{\partial v_c}{\partial y_1}&=\frac{1}{2}(x_3-x_2), \qquad \frac{\partial v_c}{\partial y_2}&=\frac{1}{2}(x_1-x_3), \qquad \frac{\partial v_c}{\partial y_3}&=\frac{1}{2}(x_2-x_1).
	\end{alignat*}
	Thus, we obtain
	\begin{alignat*}{2}
	v_c^{n+1}-v_c^n&=\frac{x_1^{n+1}-x_1^n}{\Delta t}\int_{t^n}^{t^{n+1}}\frac{1}{2}(y_2-y_3)\,\mathrm{d}t&+\frac{x_2^{n+1}-x_2^n}{\Delta t}\int_{t^n}^{t^{n+1}}\frac{1}{2}(y_3-y_1)\,\mathrm{d}t\\
	&+\frac{x_3^{n+1}-x_3^n}{\Delta t}\int_{t^n}^{t^{n+1}}\frac{1}{2}(y_1-y_2)\,\mathrm{d}t&+\frac{y_1^{n+1}-y_1^n}{\Delta t}\int_{t^n}^{t^{n+1}}\frac{1}{2}(x_3-x_2)\,\mathrm{d}t\\&+\frac{y_2^{n+1}-y_2^n}{\Delta t}\int_{t^n}^{t^{n+1}}\frac{1}{2}(x_1-x_3)\,\mathrm{d}t&+\frac{y_3^{n+1}-y_3^n}{\Delta t}\int_{t^n}^{t^{n+1}}\frac{1}{2}(x_2-x_1)\,\mathrm{d}t.
	\end{alignat*}
	
	By comparing the coefficients of $x_p^{n+1}-x_p^n$ and $y_p^{n+1}-y_p^n$ in the expressions for the algebraic and geometric volume changes, we need
	\begin{align*}
	\int_{t^n}^{t^{n+1}}\frac{1}{2}(x_p-x_q)\,\mathrm{d}t=\frac{1}{2}\Delta t(x_p-x_q), \qquad \int_{t^n}^{t^{n+1}}\frac{1}{2}(y_p-y_q)\,\mathrm{d}t=\frac{1}{2}\Delta t(y_p-y_q)
	\end{align*}
	to achieve GCL consistency. However, we have
	\begin{align}
	\begin{split}
	\int_{t^n}^{t^{n+1}}\frac{1}{2}(x_p-x_q)\,\mathrm{d}t&=\frac{1}{2}\Delta t(x_p-x_q)+\frac{1}{4}(\Delta t)^2(u_p^x-u_q^x)\\&=\frac{1}{2}\Delta t(x_p-x_q)+\frac{1}{4}\Delta t(x_p^{n+1}-x_p^n-x_q^{n+1}+x_q^n), \\
	\int_{t^n}^{t^{n+1}}\frac{1}{2}(y_p-y_q)\,\mathrm{d}t&=\frac{1}{2}\Delta t(y_p-y_q)+\frac{1}{4}(\Delta t)^2(u_p^y-u_q^y)\\&=\frac{1}{2}\Delta t(y_p-y_q)+\frac{1}{4}\Delta t(y_p^{n+1}-y_p^n-y_q^{n+1}+y_q^n).
	\end{split}
	\label{eq:agv}
	\end{align}
	GCL consistency is achieved only if the (sum of the) additional terms in \eqref{eq:agv} vanish. Specifically, for first-order time integration, GCL consistency is ensured only for simple translations of the triangle. The additional terms in \eqref{eq:agv} represent the consistency error of the forward Euler method. Therefore, this analysis confirms that second-order time accuracy is necessary to achieve GCL consistency.
	
	\begin{theorem*}
		The GCL analysis for triangles can be directly extended to quadrilaterals with straight edges. For such quadrilaterals, the geometric volume is given by
		\begin{align*}
			v_c=\frac{1}{2}\det \begin{pmatrix}
				1 & 0 & x_1 & y_1 \\ 1 & 1 & x_2 & y_2 \\ 1 & 0 & x_3 & y_3 \\ 1 & 1 & x_4 & y_4
			\end{pmatrix}
			=\frac{1}{2}(x_1y_2-x_1y_4-x_2y_1+x_2y_3-x_3y_2+x_3y_4+x_4y_1-x_4y_3).
		\end{align*}
		Comparing this geometric volume $v_c$ with the algebraic volume $\tilde{v}_c$ computed by \eqref{eq:algvol} demonstrates once again the need for using a second-order accurate time integrator to maintain GCL consistency.
	\end{theorem*}
    
	\bibliographystyle{plain}
	\bibliography{paper_lanl_bib}
    
\end{document}